\documentclass[10pt,reqno]{amsart}
\usepackage{amsmath}
\usepackage{amssymb}
\usepackage{amsthm}
\usepackage{graphicx}
\newtheorem{thm}{Theorem}[section]
\newtheorem{coro}[thm]{Corollary}
\newtheorem{lem}[thm]{Lemma}
\newtheorem{prop}[thm]{Proposition}
\theoremstyle{definition}

\theoremstyle{remark}
\newtheorem{rk}[thm]{Remark}

\numberwithin{equation}{section}

\begin{document}

\title[Inhomogeneous Strichartz estimates]
{On inhomogeneous Strichartz estimates for fractional Schr\"odinger equations
and their applications}

\author{Chu-Hee Cho, Youngwoo Koh and Ihyeok Seo}

\thanks{The first author was supported by POSTECH Math BK21 Plus Organization.} 
\thanks{The second author was supported by NRF grant 2012-008373.}

\address{Department of Mathematics, Pohang University of Science and Technology, Pohang 790-784,
Republic of Korea}
\email{chcho@postech.ac.kr}

\address{School of Mathematics, Korea Institute for Advanced Study, Seoul 130-722, Republic of Korea}
\email{ywkoh@kias.re.kr}

\address{Department of Mathematics, Sungkyunkwan University, Suwon 440-746, Republic of Korea}
\email{ihseo@skku.edu}

\maketitle

\begin{abstract}
In this paper we obtain some new inhomogeneous Strichartz estimates
for the fractional Schr\"odinger equation in the radial case.
Then we apply them to the well-posedness theory for the equation
$i\partial_{t}u+|\nabla|^{\alpha}u=V(x,t)u$, $1<\alpha<2$,
with radial $\dot{H}^\gamma$ initial data below $L^2$
and radial potentials $V\in L_t^rL_x^w$ under the scaling-critical range $\alpha/r+n/w=\alpha$.
\end{abstract}


\section{Introduction}

To begin with, let us consider the following Cauchy problem
    \begin{equation}\label{potential}
    \begin{cases}
    i\partial_{t}u + |\nabla|^{\alpha} u = F(x,t),\quad1<\alpha<2,\\
    u(x,0)=f(x),
    \end{cases}
    \end{equation}
associated with the fractional Schr\"odinger equation
 \begin{equation}\label{fse}
 i\partial_{t}u + |\nabla|^{\alpha} u = V(x,t)u
 \end{equation}
where $V:\mathbb{R}^{n+1}\rightarrow\mathbb{C}$ is a potential.
This equation has recently attracted interest from mathematical physics.
This is because fractional quantum mechanics introduced by Laskin \cite{L} is
governed by the equation where it is conjectured that physical realizations may be limited to the cases of $1<\alpha<2$.
Of course, the case $\alpha=2$ corresponds to the ordinary quantum mechanics.

By Duhamel's principle, the solution of \eqref{potential} is given by
    \begin{equation}\label{Duhamel}
    u(x,t)= e^{it |\nabla|^\alpha} f(x) - i \int_{0}^{t} e^{i(t-s) |\nabla|^\alpha} F(\cdot,s)ds,
    \end{equation}
where the propagator $e^{it |\nabla|^\alpha}$ is given by means of the Fourier transform, as follows:
$$e^{it |\nabla|^\alpha} f(x)= \frac{1}{(2\pi)^n} \int_{\mathbb{R}^n} e^{ix\cdot\xi + it|\xi|^\alpha}\widehat{f}(\xi)d\xi.$$
Then the standard approach to the problem \eqref{potential}
is to obtain the corresponding Strichartz estimates
which control space-time integrability of \eqref{Duhamel}
in view of that of the initial datum $f$ and the forcing term $F$.

In the classical case $\alpha=2$, the Strichartz estimates originated by Strichartz \cite{Str} have been extensively studied by many authors (\cite{GV,KT,CW,K,F,V,Ko,LS,LS2,CN,CN2,S}).
Over the past several years, considerable attention has been paid to the fractional order where $1<\alpha<2$
in the radial case (see \cite{Sh,GW,Ke} and references therein).
From these works, when $2n/(2n-1)\leq\alpha<2$, the homogeneous Strichartz estimate
 \begin{equation}\label{homo}
    \|e^{it |\nabla|^\alpha} f\|_{L^{q}_{t} L^{p}_{x} }
    \lesssim \|f\|_{\dot{H}^\gamma}
    \end{equation}
holds for radial functions $f\in\dot{H}^\gamma(\mathbb{R}^n)$ if
\begin{equation}\label{con}
\frac{\alpha}{q}+\frac{n}{p}=\frac{n}{2}-\gamma,\quad 2 \leq q \leq \infty \quad\mbox{and}\quad (q,p) \neq (2, \frac{4n-2}{2n-3}).
\end{equation}
Here the condition \eqref{con} is optimal if $2n/(2n-1)<\alpha<2$.
But when $\alpha=2n/(2n-1)$, \eqref{homo} is unknown for the endpoint $(q,p)=(2,(4n-2)/(2n-3))$.
Also, it is known that the estimate does not hold in general if $f$ does not have radial symmetry.

\begin{figure}[t!]\label{figure}
\includegraphics[height=6cm]{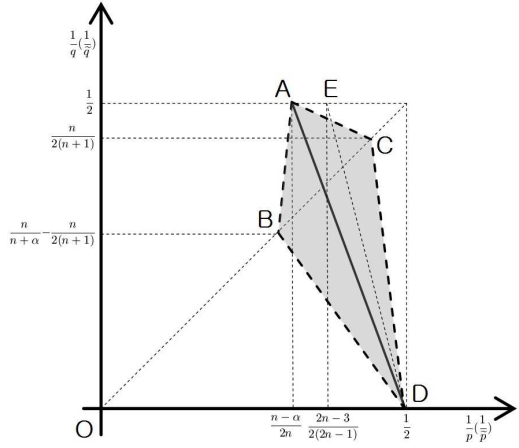}
\caption{The range for Corollary \ref{mixedprop}.
Here the open segment $(B,C)$ is the range for \eqref{inhomo2},
and $[A,D]$ is the range for \eqref{inhomo}.}\label{fig1}
\end{figure}

Now, by duality and the Christ-Kiselev lemma (\cite{CK}),
one may use \eqref{homo} with $\gamma=0$ to get some inhomogeneous estimates
 \begin{equation}\label{inhomo}
    \Big\| \int_{0}^{t}e^{i(t-s) |\nabla|^\alpha} F(\cdot,s) ds \Big\|_{L^{q}_t L^{p}_x}
    \lesssim \|F\|_{L^{\widetilde{q}'}_tL^{\widetilde{p}'}_x}
    \end{equation}
for $(q,p)$ and $(\widetilde{q},\widetilde{p})$ which satisfy \eqref{con} with $\gamma=0$ and $q>\widetilde{q}'$.
This means that $(q,p)$ and $(\widetilde{q},\widetilde{p})$ lie on the segment $AD$ in Figure \ref{figure}.
However, these trivial estimates are not enough to imply the well-posedness for the equation \eqref{fse}
with the initial data $f\in\dot{H}^\gamma(\mathbb{R}^n)$ beyond the case $\gamma=0$.
When $\gamma\neq0$ we need to obtain \eqref{inhomo} on a wider range of $(q,p)$ and $(\widetilde{q},\widetilde{p})$.
See Section \ref{sec2} for details.

Let us first mention the following necessary conditions for \eqref{inhomo}:
    \begin{equation}\label{nece}
    \alpha(\frac{1}{q}+\frac{1}{\widetilde{q}})- n(1-\frac{1}{p}-\frac{1}{\widetilde{p}}) =0
    \end{equation}
and
    \begin{equation}\label{nece2}
    \frac{1}{q}+\frac{n}{p} < \frac{n}{2},\quad\frac{1}{\widetilde{q}}+\frac{n}{\widetilde{p}}< \frac{n}{2}.
    \end{equation}
The first is just the scaling condition and the second will be shown in Section \ref{sec4}.

Our main result in this paper is the following theorem where we obtain \eqref{inhomo}
on the open segment $BC$ in Figure \ref{figure}.

\begin{thm}\label{prop}
Let $n\geq2$ and $2n/(2n-1)\leq\alpha<2$.
Assume that $F(x,t)$ is a radial function with respect to the spatial variable $x$.
Then we have
    \begin{equation}\label{inhomo2}
    \Big\| \int_{0}^{t} e^{i(t-s) |\nabla|^\alpha} F(\cdot,s) ds \Big\|_{L^{q}_{x,t} }
    \lesssim \|F\|_{L^{\widetilde{q}'}_{x,t} }
    \end{equation}
if
    \begin{equation}\label{sharp}
    \frac{1}{q},\frac{1}{\widetilde{q}} < \frac{n}{2(n+1)}
    \quad\mbox{and}\quad \frac{1}{q}+\frac{1}{\widetilde{q}} = \frac{n}{n+\alpha}.
    \end{equation}
\end{thm}

\begin{rk}\label{rem}
It should be noted that the range \eqref{sharp} is sharp.
Namely, the second condition in \eqref{sharp} is the scaling condition for \eqref{inhomo2} (see \eqref{nece}),
and the first one is the necessary condition \eqref{nece2} when $q=p$ and $\widetilde{q}=\widetilde{p}$.
\end{rk}


From interpolation between \eqref{inhomo} and \eqref{inhomo2},
we can directly obtain further estimates when $(q,p)$ and $(\widetilde{q},\widetilde{p})$
are contained in the open quadrangle with vertices $A,B,D,C$.
Precisely, we have the following corollary.

\begin{coro}\label{mixedprop}
Let $n\geq2$ and $2n/(2n-1)\leq\alpha<2$.
Assume that $F(x,t)$ is a radial function with respect to the spatial variable $x$,
and that $(q,p)$ and $(\widetilde{q},\widetilde{p})$ satisfy
the necessary conditions \eqref{nece} and \eqref{nece2}.
Then we have
    \begin{equation}\label{inhs}
    \Big\| \int_{0}^{t}e^{i(t-s) |\nabla|^\alpha} F(\cdot,s) ds \Big\|_{L^{q}_t L^{p}_x}
    \lesssim \|F\|_{L^{\widetilde{q}'}_tL^{\widetilde{p}'}_x}
    \end{equation}
if the following conditions hold:
\begin{itemize}
\item For $(q,p)$,
    \begin{equation}\label{freq_cond_2}
    \frac{-n(n+2-\alpha)}{(2\alpha-1)n+\alpha} (\frac{1}{p}-\frac{1}{2})    <\frac{1}{q} < \frac{-n}{(\alpha-1)n+\alpha} (\frac{1}{p}-\frac{n-\alpha}{2n}) +\frac{1}{2},
    \end{equation}
and
    \begin{equation}\label{freq_cond_4}
    \frac{1}{p} > \frac{(\alpha-1)n^2 -\alpha^2 n -\alpha^2}{((2\alpha-1)n+\alpha)n}\frac{1}{q}
    +\frac{\alpha(n+2-\alpha)}{2((2\alpha-1)n +\alpha)} .
    \end{equation}
\item Similarly for $(\widetilde{q},\widetilde{p})$.
\end{itemize}
\end{coro}

Let us give more details about the conditions in the above corollary.
The line $BD$ in Figure \ref{fig1} is when the equality holds in the first inequality of \eqref{freq_cond_2}.
Similarly, the lines $AC$ and $AB$ correspond to the second inequality in \eqref{freq_cond_2} and the inequality \eqref{freq_cond_4},
respectively.
Finally, the line $CD$ is sharp because it is determined from the necessary condition \eqref{nece2}.

\medskip

Now we apply the above Strichartz estimates to the well-posedness theory
for the fractional Schr\"odinger equation in the radial case:
\begin{equation}\label{cauchy}
\begin{cases}
i\partial_{t}u+|\nabla|^{\alpha}u=V(x,t)u,\quad1<\alpha<2,\\
u(x,0)=f(x),
\end{cases}
\end{equation}
where we assume that $u$, $f$ and $V$ are radial functions with respect to the spatial variable $x$.
We obtain the following well-posedness for \eqref{cauchy}
with $\dot{H}^\gamma$ initial data $f$ below $L^2$
and potentials $V\in L_t^rL_x^w$ under the scaling-critical range $\alpha/r+n/w=\alpha$.
The Cauchy problem \eqref{cauchy} was studied in \cite{D'APV,NS} particularly when $\alpha=2$ and $\gamma=0$.

\begin{thm}\label{thm}
Let $\frac{2n}{2n-1}\leq\alpha<2$ and $\frac{-(\alpha-1)n}{2(n+1)}<\gamma\leq0$ for $n\geq2$.
Assume that $f\in\dot{H}^\gamma(\mathbb{R}^n)$ and
$V\in L_t^{r}([0,T];L_x^w(\mathbb{R}^n))$ for some $T>0$.
Then there exists a unique solution
$u\in C([0,T];\dot{H}^\gamma(\mathbb{R}^n))\cap L_t^q([0,T];L_x^p(\mathbb{R}^n))$ of \eqref{cauchy}
if
\begin{equation}\label{cond}
\frac{\alpha}{q}+\frac{n}{p}=\frac{n}2-\gamma,\quad\frac{\alpha}r+\frac{n}w=\alpha,
\end{equation}
\begin{equation}\label{condcon}
\frac{-\gamma}{\alpha-1}<\frac1q<\frac{\gamma}{(\alpha-1)n}+\frac12,
\end{equation}
and
\begin{equation}\label{cond4}
1-\frac1q-\frac{n(n+\alpha)(\alpha-1)+2\gamma((2\alpha-1)n+\alpha)}{2n(n+\alpha)(\alpha-1)}
<\frac1r<1-\frac1q+\frac{\gamma(n+2-\alpha)}{(n+\alpha)(\alpha-1)}.
\end{equation}
\end{thm}

\begin{rk}
The condition $\alpha/r+n/w=\alpha$ on the potential is critical in the sense of scaling.
Indeed, $u_{\epsilon}(x,t)=u(\epsilon x, \epsilon^\alpha t)$ takes \eqref{cauchy} into  $i\partial_{t}u_{\epsilon}+|\nabla|^{\alpha}u_{\epsilon}=V_{\epsilon}(x,t)u_{\epsilon}$ with $V_{\epsilon}(x,t)=\epsilon^{\alpha}V(\epsilon x,\epsilon^{\alpha}t)$.
Hence the norm
$$\|V_{\epsilon}\|_{L_t^{r}L_x^w}=
\|V_{\epsilon}\|_{L_t^{r}(\mathbb{R};L_x^w(\mathbb{R}^n))}
=\epsilon^{\alpha-\alpha/r-n/w}\|V\|_{L_t^{r}L_x^w}$$
is independent of $\epsilon$ precisely when $\alpha/r+n/w=\alpha$.
\end{rk}

\begin{rk}
In Proposition 3.9 of \cite{GW},
the inhomogeneous estimates were shown in certain region that lies below the segment $ED$ in Figure \ref{figure}.
(Note that $(1/q,1/p)\in ED$ if and only if $q,p\geq2$ and $2/q+(2n-1)/p=n-1/2$.
Also, the point $E$ is the same as $A$ when $\alpha=2n/(2n-1)$.)
By interpolation between these and our estimates, we can also obtain further estimates
in the triangle with vertices $A,C,E$.
We omit the details since it does not affect the range $\frac{-(\alpha-1)n}{2(n+1)}<\gamma\leq0$
of $\gamma$ in Theorem \ref{thm} (see Section \ref{sec2}).
\end{rk}

The rest of the paper is organized as follows.
In Section \ref{sec2}, we prove Theorem \ref{thm}
by making use of the Strichartz estimates \eqref{homo} and \eqref{inhs}.
Section \ref{sec3} is devoted to proving Theorem \ref{prop},
and in Section \ref{sec4} we show the necessary condition \eqref{nece2}.
In the final section, Section \ref{sec5},
we show Lemma \ref{bessel} which gives some estimates for Bessel functions
and is used for the proof of Theorem \ref{prop}.

\medskip

Throughout the paper, we shall use the letter $C$ to denote positive constants
which may be different at each occurrence.
We also use the symbol $\widehat{f}$ to denote the Fourier transform of $f$,
and denote $A\lesssim B$ and $A\sim B$ to mean $A\leq CB$ and $CB\leq A\leq CB$, respectively,
with unspecified constants $C>0$.

\section{Application}\label{sec2}
In this section we prove Theorem \ref{thm}.
The proof is quite standard but we need to observe that
if $(q,p)$ and $(\widetilde{q},\widetilde{p})$ satisfy the inhomogeneous estimate \eqref{inhs},
then the midpoint of them lies on the segment $AD$.
Note that $(1/q,1/p)\in AD$ if and only if $q,p\geq2$ and $\alpha/q+n/p=n/2$.
Hence, if $\alpha/q+n/p=n/2-\gamma$ for $\gamma\in\mathbb{R}$,
then $(\widetilde{q},\widetilde{p})$ should satisfy
$\alpha/\widetilde{q}+n/\widetilde{p}=n/2+\gamma$ to give \eqref{inhs}.
In what follows, it will be convenient to keep in mind this key observation.

By Duhamel's principle, the solution of \eqref{cauchy} is given by
\begin{equation}\label{phi}
\Phi(u):=e^{it|\nabla|^\alpha}f(x)-i\int_{0}^{t} e^{i(t-s)|\nabla|^\alpha}V(\cdot,s)u(\cdot,s)ds.
\end{equation}
Then the standard fixed-point argument is to choose the solution space on which $\Phi$ is a contraction mapping. The Strichartz estimates play a central role in this step.
Indeed, by the estimates \eqref{homo} and \eqref{inhs}, we see that
\begin{equation}\label{ma}
\|\Phi(u)\|_{L_t^q([0,T];L_x^p)}\leq C\|f\|_{\dot{H}^\gamma}+C\|Vu\|_{L^{\widetilde{q}'}_t([0,T]; L^{\widetilde{p}'}_x)}
\end{equation}
if
\begin{equation}\label{cond5}
\frac{-(\alpha-1)n}{2(n+1)}<\gamma\leq0
\end{equation}
\begin{equation}\label{cond6}
\frac{\alpha}{q}+\frac{n}{p}=\frac{n}2-\gamma,
\quad\frac{\alpha}{\widetilde{q}}+\frac{n}{\widetilde{p}}=\frac{n}2+\gamma,
\end{equation}
\begin{equation}\label{cond7}
\frac{-\gamma}{\alpha-1}<\frac1q<\frac{\gamma}{(\alpha-1)n}+\frac12,
\end{equation}
and
\begin{equation}\label{cond8}
\frac{\gamma(n+2-\alpha)}{(n+\alpha)(\alpha-1)}
<\frac1{\widetilde{q}}
<\frac{n(n+\alpha)(\alpha-1)+2\gamma((2\alpha-1)n+\alpha)}{2n(n+\alpha)(\alpha-1)}.
\end{equation}
Here, the conditions \eqref{cond5} and \eqref{cond7} are given from that
the line $\frac{\alpha}{q}+\frac{n}{p}=\frac{n}2-\gamma$ lies in the closed triangle
with vertices $A,C,D$ except the closed segments $[A,C],[C,D]$.
Note that $\gamma=\frac{-(\alpha-1)n}{2(n+1)}$
when this line passes through the point $C$.
Similarly, the condition \eqref{cond8} is given from that
the line $\frac{\alpha}{\widetilde{q}}+\frac{n}{\widetilde{p}}=\frac{n}2+\gamma$ lies in the closed triangle
with vertices $A,B,D$ except the closed segments $[A,B],[B,D]$.

By H\"older's inequality, we now get
$$\|\Phi(u)\|_{L_t^q([0,T];L_x^p)}
\leq C\|f\|_{\dot{H}^\gamma}+C\|V\|_{L^r_t([0,T];L^w_x)}
\|u\|_{L^q_t([0,T];L^p_x)}$$
if $\alpha/r+n/w=\alpha$ and the condition \eqref{cond4} holds.
Indeed, when applying H\"older's inequality to the second term in the right-hand side of \eqref{ma},
the conditions $\alpha/r+n/w=\alpha$ and \eqref{cond4}
follow from \eqref{cond6} and \eqref{cond8}, respectively.

From the above argument and the linearity, it follows that
\begin{align*}
\|\Phi(u)-\Phi(v)\|_{L_t^q([0,T];L_x^p)}
&\leq C\|V\|_{L^r_t([0,T];L^w_x)}\|u-v\|_{L^q_t([0,T];L^p_x)}\\
&\leq\frac12\|u-v\|_{L^q_t([0,T];L^p_x)},
\end{align*}
which says that $\Phi$ is a contraction mapping,
if $T$ is sufficiently small.
But here, since the above process works also on time-translated small intervals
if $u(\cdot,t)\in\dot{H}^\gamma(\mathbb{R}^n)$ for all $t\geq0$,
the smallness assumption on $T$ can be removed by
iterating the process a finite number of times.
For this we will show that
\begin{equation}\label{final}
\|u\|_{L_t^\infty\dot{H}_x^\gamma}\lesssim\|f\|_{\dot{H}^\gamma}
+\|V\|_{L^r_tL^w_x}\|u\|_{L^q_tL^p_x}.
\end{equation}
From \eqref{phi}, we first see that
$$\|u\|_{L_t^\infty\dot{H}_x^\gamma}\lesssim\|e^{it|\nabla|^\alpha}D^\gamma f\|_{L_t^\infty L_x^2}
+\Big\|\int_{0}^{t} e^{i(t-s)|\nabla|^\alpha}D^\gamma(V(\cdot,s)u(\cdot,s))ds\Big\|_{L_t^\infty L_x^2}.$$
Since $e^{it|\nabla|^\alpha}$ is an isometry in $L^2$,
the first term in the right-hand side is clearly bounded by $\|f\|_{\dot{H}^\gamma}$.
On the other hand, by the inhomogeneous estimate \eqref{inhomo} the second term is bounded by
$\|D^\gamma(Vu)\|_{L_t^{\widetilde{u}'}L_x^{\widetilde{v}'}}$,
where $\widetilde{u},\widetilde{v}\geq2$ and $\alpha/\widetilde{u}+n/\widetilde{v}=n/2$.
Here we use the Sobolev embedding
$$\|g\|_{L^b}\lesssim\|D^\beta g\|_{L^a},$$
where $1/a-1/b=\beta/n$ with $0\leq\beta<n/a$ and $1<a<\infty$,
and H\"older's inequality to get
\begin{align*}
\|D^\gamma(Vu)\|_{L_t^{\widetilde{u}'}L_x^{\widetilde{v}'}}
&\lesssim\|Vu\|_{L_t^{\widetilde{u}'}L_x^{a}}\\
&\leq\|V\|_{L_t^{r}L_x^{w}}\|u\|_{L_t^{q}L_x^{p}}.
\end{align*}
The required conditions here are summarized as follows:
$$\widetilde{u},\widetilde{v}\geq2,\quad\frac{\alpha}{\widetilde{u}}+\frac n{\widetilde{v}}=\frac n2,$$
$$\frac1a-\frac1{\widetilde{v}'}=\frac{-\gamma}n,\quad0\leq-\gamma<\frac na,\quad1<a<\infty,$$
$$\frac1{\widetilde{u}'}=\frac1r+\frac1q,\quad\frac1a=\frac1w+\frac1p.$$
But, the inequalities $\widetilde{u},\widetilde{v}\geq2$ and $1<a<\infty$ are satisfied automatically
from the conditions on $q,r,p,w$ in Theorem \ref{thm}.
On the other hand, the inequality $0\leq-\gamma<\frac na$ is redundant
because $\widetilde{v}\geq2$.
The remaining four equalities is reduced to the following one equality
$$\frac{\alpha}{q}+\frac{n}{p}-\frac{n}2+\gamma=-(\frac{\alpha}{r}+\frac{n}{w}-\alpha)$$
which is clearly satisfied from the condition \eqref{cond}.
Consequently, we get \eqref{final}.

\section{Inhomogeneous Strichartz estimates}\label{sec3}

In this section we prove Theorem \ref{prop}.
Let us first consider the multiplier operators $P_k$ for $k \in \mathbb{Z}$
defined by
\begin{equation*}
\widehat{P_kf} = \phi(|\cdot|/2^k) \widehat{f},
\end{equation*}
where $\phi:\mathbb{R}\rightarrow[0,1]$ is a smooth cut-off function
which is supported in $(1/2,2)$ and
satisfies $\sum_{k \in \mathbb{Z}} \phi (\cdot/2^k) =1$.
Then we will obtain the following frequency localized estimates
(Proposition \ref{prop2})
which imply Theorem \ref{prop}.

\begin{prop}\label{prop2}
Let $n\geq2$ and $2n/(2n-1)\leq\alpha<2$.
Assume that $F(x,t)$ is a radial function with respect to the spatial variable $x$.
Then we have
    \begin{equation}\label{inhomo3}
    \Big\| \int_{\mathbb{R}} e^{i(t-s) |\nabla|^\alpha} P_k F(\cdot,s) ds \Big\|_{L^{q}_{x,t} }
    \lesssim \|F\|_{L^{\widetilde{q}'}_{x,t} }
    \end{equation}
uniformly in $k \in \mathbb{Z}$ if
    \begin{equation}\label{accep_and_scal_cond}
     \frac{1}{q}, \frac{1}{\widetilde{q}} < \frac{n}{2(n+1)}
    \quad\mbox{and}\quad \frac{1}{q} + \frac{1}{\widetilde{q}} = \frac{n}{n+\alpha}.
    \end{equation}
\end{prop}

Indeed, since $q >2$ from the first condition in \eqref{sharp},
by the Littlewood-Paley theorem and Minkowski integral inequality,
one can see that
    \begin{align*}
    \Big\| \int_{\mathbb{R}} e^{i(t-s) |\nabla|^\alpha} F(\cdot,s) ds \Big\|_{L^{q}_{x,t} }^2
    &\leq C\Big\| \Big(\sum_{k\in\mathbb{Z}}\Big|\int_{\mathbb{R}} e^{i(t-s)|\nabla|^\alpha} P_kF(\cdot,s) ds\Big|^2\Big)^{1/2} \Big\|_{L^{q}_{x,t} }^2\\
    &\leq C\sum_{k\in\mathbb{Z}}\Big\| \int_{\mathbb{R}} e^{i(t-s) |\nabla|^\alpha} P_k\big(\sum_{|j-k|\leq1}P_j F(\cdot,s)\big)ds \Big\|_{L^{q}_{x,t} }^2.
    \end{align*}
Now, by \eqref{inhomo3}, the right-hand side in the above is bounded by
    $$
    C\sum_{k\in\mathbb{Z}}\big\|\sum_{|j-k|\leq1}P_j F\big\|_{L^{\widetilde{q}'}_{x,t}}^2.
    $$
Since $\widetilde{q}'< 2$,
using the Minkowski integral inequality and Littlewood-Paley theorem,
this is bounded by $C\|F\|_{L^{\widetilde{q}'}_{x,t}}^2$.
By this boundedness and $\widetilde{q}'< 2<q$, one may now use
the Christ-Kiselev lemma (\cite{CK}) to obtain
    $$
    \Big\| \int_{0}^{t}e^{i(t-s) |\nabla|^\alpha} F(\cdot,s) ds \Big\|_{L^{q}_{x,t}}
    \lesssim \|F\|_{L^{\widetilde{q}'}_{x,t}}
    $$
as desired.

Now it remains to prove the above proposition.

\subsection{Proof of Proposition \ref{prop2}}

Since we are assuming the scaling condition in \eqref{accep_and_scal_cond},
by rescaling $(x,t)\rightarrow(\lambda x,\lambda^{\alpha}t)$,
we may show \eqref{inhomo3} only for $k=0$.

Let us first consider
$x= rx'$, $y=\lambda y'$ and $\xi=\rho\xi'$ for $x', y', \xi' \in S^{n-1}$,
where $r = |x|$, $\lambda=|y|$ and $\rho = |\xi|$.
Then by using the fact (see \cite{St}, p. 347) that
    $$
    \int_{S^{n-1}} e^{-ir\rho x'\cdot\xi'} dx' = c_n (r\rho)^{-\frac{n-2}{2}} J_{\frac{n-2}{2}} (r\rho),
    $$
where $J_m$ denotes the Bessel function with order $m$,
it is easy to see that
    \begin{align}\label{be}
    &\int_{\mathbb{R}} e^{i(t-s) |\nabla|^\alpha} P_0 F(\cdot,s) ds \\
    \nonumber&=\int_{\mathbb{R}} \int_{\mathbb{R}} \int_{\mathbb{R}} \int_{S^{n-1}} \int_{S^{n-1}}
        e^{i\rho(rx'-\lambda y')\cdot\xi' + i(t-s)\rho^\alpha} \phi(\rho) F(\lambda y',s)
        (\lambda\rho)^{n-1} dy' d\xi' d\rho d\lambda ds\\
    \nonumber&=r^{-\frac{n-2}{2}} \int_{\mathbb{R}} \lambda^{-\frac{n-2}{2}}
        \bigg( \int_{\mathbb{R}} e^{i(t-s)\rho^\alpha} J_{\frac{n-2}{2}} (r\rho) J_{\frac{n-2}{2}} (\lambda\rho) \rho \phi(\rho) d\rho \bigg)
        \big[\lambda^{n-1}F(\lambda y',s) \big] d\lambda ds .
    \end{align}
Now we define the operators $T_j h$, $j\geq0$, as
    \begin{equation}\label{110}
    T_0 h (r,t)
    = \chi_{(0,1)}(r)~ r^{-\frac{n-2}{2}} \int_0^\infty e^{it\rho^\alpha} J_{\frac{n-2}{2}} (r\rho) \varphi(\rho) h(\rho) d\rho
    \end{equation}
and for $j \geq1$
   \begin{equation}\label{111}
    T_j h (r,t)
    = \chi_{[2^{j-1},2^j)}(r)~ r^{-\frac{n-2}{2}} \int_0^\infty e^{it\rho^\alpha} J_{\frac{n-2}{2}} (r\rho) \varphi(\rho) h(\rho) d\rho,
    \end{equation}
where $\chi_A$ denotes the characteristic function of a set $A$ and
$\varphi^2(\rho) = \rho \phi(\rho)$.
Then the adjoint operator $T^*_{k}$ of $T_{k}$ is given by
     $$
    T^{*}_{0} H (\rho)
    = \varphi (\rho) \int_{\mathbb{R}}e^{-is\rho^{\alpha}}\int_0^\infty\chi_{(0,1)}(\lambda)~ \lambda^{-\frac{n-2}{2}} J_{\frac{n-2}{2}} (\lambda\rho) H(\lambda,s) d\lambda ds
    $$
and for $k\geq1$
    \begin{equation}\label{112}
    T^{*}_{k} H (\rho)
    = \varphi (\rho) \int_{\mathbb{R}}e^{-is\rho^{\alpha}}\int_0^\infty\chi_{[2^{k-1},2^k)}(\lambda)~ \lambda^{-\frac{n-2}{2}} J_{\frac{n-2}{2}} (\lambda\rho) H(\lambda,s) d\lambda ds.
    \end{equation}
Hence, since $F(\lambda y',s)$ is independent of $y'\in S^{n-1}$,
by setting $H(\lambda,s)=F(\lambda y',s)$, it follows from \eqref{be} that
\begin{equation*}
    \int_{\mathbb{R}} e^{i(t-s) |\nabla|^\alpha} P_0 F(\cdot,s) ds = \sum_{j,k \geq0} T_j T^{*}_{k} (\lambda^{n-1}H).
    \end{equation*}
Now we are reduced to showing that
    \begin{equation}\label{inhomo_est_prop}
    \Big\| \sum_{j,k \geq 0} T_j T^{*}_{k} (\lambda^{n-1}H) \Big\|_{L^{q}_{t} \mathfrak{L}^{q}_{r}}
    \lesssim \|H\|_{L^{\widetilde{q}'}_t \mathfrak{L}^{\widetilde{q}'}_r},
    \end{equation}
where we denote by $\mathfrak{L}^q_r$ the space $L^q(r^{n-1}dr)$.

\begin{figure}[t!]
\includegraphics[height=6cm]{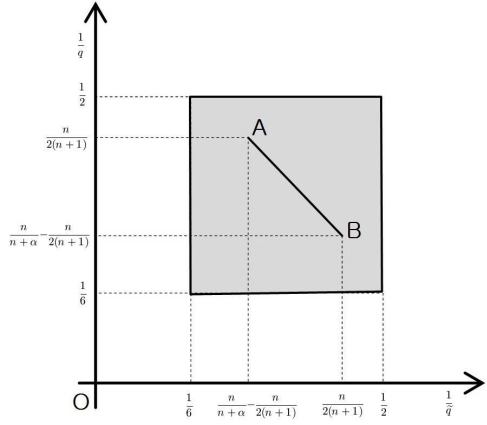}
\caption{The range of $q,\widetilde{q}$ for Lemma \ref{local_prop}.
Here the open segment $(A,B)$ is the range for Proposition \ref{prop2} (see \eqref{accep_and_scal_cond}).}\label{fig2}
\end{figure}

From now on, we will show \eqref{inhomo_est_prop} by making use of the following lemma
which will be obtained in Subsection \ref{subsec3.2}.

\begin{lem}\label{local_prop}
Let $n\geq2$ and $2n/(2n-1)\leq\alpha<2$.
Then we have for $j,k \geq 0$
    \begin{align*}
    \| T_j T^{*}_{k} &(\lambda^{n-1}H) \|_{L^{q}_{t} \mathfrak{L}^{q}_{r}}\\
    &\lesssim 2^{j(\frac{2n+1}{2q}-\frac{2n-1}{4})} 2^{k(\frac{2n+1}{2\widetilde{q}} - \frac{2n-1}{4})} 2^{\frac{-|j-k|}{2}(\frac{1}{2}-\max(\frac{1}{q},\frac{1}{\widetilde{q}}))}
    \|H\|_{L^{\widetilde{q}'}_t \mathfrak{L}^{\widetilde{q}'}_r}
    \end{align*}
if $2\leq q,\widetilde{q}\leq6$ (see Figure \ref{fig2}).
\end{lem}

\noindent\textbf{The case $2n/(2n-1) <\alpha<2$.}
We first decompose the sum over $j,k$ into two parts, $j\le k$ and $j\ge k$:
    \begin{align*}
\sum_{j, k\geq 0} &\| T_j T^{*}_{k} (\lambda^{n-1}H) \|_{L^{q}_{t} \mathfrak{L}^{q}_{r}} \\
    &\leq \sum_{j=0}^{\infty} \sum_{k=j}^{\infty} \| T_j T^{*}_{k} (\lambda^{n-1}H) \|_{L^{q}_{t} \mathfrak{L}^{q}_{r}}
        + \sum_{k=0}^{\infty} \sum_{j=k}^{\infty} \| T_j T^{*}_{k} (\lambda^{n-1}H) \|_{L^{q}_{t} \mathfrak{L}^{q}_{r}} .
    \end{align*}
When $j\le k$, using Lemma \ref{local_prop}, we then have
    \begin{align*}
\sum_{j=0}^{\infty} \sum_{k=j}^{\infty} &\| T_j T^{*}_{k} (\lambda^{n-1}H) \|_{L^{q}_{t} \mathfrak{L}^{q}_{r}} \\
    &\lesssim \sum_{j=0}^{\infty}\sum_{k=j}^{\infty}
        2^{j(\frac{2n+1}{2q}-\frac{2n-1}{4})}
        2^{k(\frac{2n+1}{2\widetilde{q}}-\frac{2n-1}{4})}
        2^{\frac{-|j-k|}{2}(\frac{1}{2}-\max(\frac{1}{q},\frac{1}{\widetilde{q}}))} \|H\|_{L^{\widetilde{q}'}_t\mathfrak{L}^{\widetilde{q}'}_r}\\
    &=\sum_{j=0}^{\infty}
        2^{j(\frac{2n+1}{2q}-\frac{2n-2}{4}-\frac{1}{2}\max(\frac{1}{q},\frac{1}{\widetilde{q}}))}
        \sum_{k=j}^{\infty}2^{k(\frac{2n+1}{2\widetilde{q}}-\frac{n}{2}+ \frac{1}{2}\max(\frac{1}{q},\frac{1}{\widetilde{q}}))}
        \|H\|_{L^{\widetilde{q}'}_t \mathfrak{L}^{\widetilde{q}'}_r} .
    \end{align*}
Note here that the first condition in \eqref{accep_and_scal_cond} implies
    \begin{equation}\label{1st_cond_modify}
    \frac{2n+1}{2\widetilde{q}} - \frac{n}{2} + \frac{1}{2}\max(\frac{1}{q},\frac{1}{\widetilde{q}})
    \leq (n+1)\max(\frac{1}{q},\frac{1}{\widetilde{q}}) - \frac{n}{2}
    <0.
    \end{equation}
From this, it follows that
\begin{align*}
\sum_{j=0}^{\infty}2^{j(\frac{2n+1}{2q}-\frac{2n-2}{4}-\frac{1}{2}
\max(\frac{1}{q},\frac{1}{\widetilde{q}}))}
\sum_{k=j}^{\infty}& 2^{k(\frac{2n+1}{2\widetilde{q}}-\frac{n}{2}+ \frac{1}{2}\max(\frac{1}{q},\frac{1}{\widetilde{q}}))}\\
&\lesssim
\sum_{j=0}^{\infty}2^{j(\frac{2n+1}{2}(\frac{1}{q}+\frac{1}{\widetilde{q}})-\frac{2n-1}{2})}.
\end{align*}
On the other hand, the second condition in \eqref{accep_and_scal_cond} implies
    \begin{equation}\label{range_B}
    \frac{2n+1}{2}(\frac{1}{q}+\frac{1}{\widetilde{q}}) - \frac{2n-1}{2} < 0
    \end{equation}
since $\alpha >2n/(2n-1)$.
Consequently, we get
\begin{align*}
    \sum_{j=0}^{\infty} \sum_{k=j}^{\infty} \| T_j T^{*}_{k} (\lambda^{n-1}H) \|_{L^{q}_{t} \mathfrak{L}^{q}_{r}}&\lesssim
\sum_{j=0}^{\infty}2^{j(\frac{2n+1}{2}(\frac{1}{q}+\frac{1}{\widetilde{q}}) - \frac{2n-1}{2})}
\|H\|_{L^{\widetilde{q}'}_t \mathfrak{L}^{\widetilde{q}'}_r}\\
    &\lesssim\|H\|_{L^{\widetilde{q}'}_t \mathfrak{L}^{\widetilde{q}'}_r}
\end{align*}
as desired.
The other part where $j\ge k$ follows clearly from the same argument.

\

\noindent\textbf{The case $\alpha = 2n/(2n-1)$.}
The previous argument is no longer available in this case,
since the left-hand side in \eqref{range_B} becomes zero.
But here we deduce \eqref{inhomo_est_prop} from bilinear interpolation
between bilinear form estimates which follow from Lemma \ref{local_prop}.
This enables us to gain some summability as before.

Let us first define the bilinear operators $B_{j,k}$ by
 \begin{equation*}
    B_{j,k} (H,\widetilde{H}) =
    \Big\langle T^*_{k} (\lambda^{n-1}H), T^*_{j} (\lambda^{n-1} \widetilde{H}) \Big\rangle_{L^2_{r,t}},
    \end{equation*}
where $\langle\,\,,\,\rangle$ denotes the usual inner product on the space $L^2_{r,t}$.
Then it is enough to show the following bilinear form estimate
    \begin{equation}\label{2:1}
    \Big| \sum_{j,k \geq 0} B_{j,k} (H,\widetilde{H}) \Big|
    \lesssim \|H\|_{L^{\widetilde{q}'}_t \mathfrak{L}^{\widetilde{q}'}_r}
    \|\widetilde{H}\|_{L^{q'}_t \mathfrak{L}^{q'}_r}.
    \end{equation}
In fact, from \eqref{2:1} we get
    \begin{align*}
    \Big\| \sum_{j,k \geq 0} T_j T^{*}_{k} (\lambda^{n-1}H) \Big\|_{L^{q}_{t} \mathfrak{L}^{q}_{r}}
    &= \sup_{\|\widetilde{H}\|_{L^{q'}_{r,t}}=1}
        \iint \sum_{j,k \geq 0} T_j T^{*}_{k} (\lambda^{n-1}H) ~ r^{\frac{n-1}{q}} \widetilde{H}(r,t) dr dt \\
    &= \sup_{\|\widetilde{H}\|_{L^{q'}_{r,t}}=1}
        \sum_{j,k \geq 0}
        \Big\langle T^*_{k} (\lambda^{n-1}H), T^*_{j} (r^{\frac{n-1}{q}} \widetilde{H}) \Big\rangle_{ L^2_{r,t}}\\
    &\lesssim \sup_{\|\widetilde{H}\|_{L^{q'}_{r,t}}=1}
        \|H\|_{L^{\widetilde{q}'}_t \mathfrak{L}^{\widetilde{q}'}_r}  \| r^{-\frac{n-1}{q'}} \widetilde{H}\|_{L^{q'}_t \mathfrak{L}^{q'}_r}
    = \|H\|_{L^{\widetilde{q}'}_t \mathfrak{L}^{\widetilde{q}'}_r}
    \end{align*}
as desired.

For \eqref{2:1}, we first decompose the sum over $j,k$ into two parts,
$j\leq k$ and $j\geq k$:
    $$
    \Big| \sum_{j,k \geq 0} B_{j,k} (H,\widetilde{H}) \Big|
    \leq \sum_{j=0}^{\infty} \Big| \sum_{k=j}^{\infty} B_{j,k} (H,\widetilde{H}) \Big|
        + \sum_{k=0}^{\infty} \Big| \sum_{j=k}^{\infty} B_{j,k} (H,\widetilde{H}) \Big|.
    $$
Then we will use the following estimate which follows from
H\"older's inequality and Lemma \ref{local_prop}:
    \begin{align*}
    |B_{j,k}(H,\widetilde{H})|
    &=\iint\widetilde{H}(r,s)r^{\frac{n-1}{q'}}r^{\frac{n-1}{q}}T_j T^{*}_{k}(\lambda^{n-1}H)drds\\
    &\leq\| \widetilde{H} \|_{L^{q'}_{t} \mathfrak{L}^{q'}_{r}}
    \big\| T_j T^{*}_{k}(\lambda^{n-1}H)\big\|_{L^{q}_{t}\mathfrak{L}^{q}_{r}}\\
    &\lesssim2^{j(\frac{2n+1}{2q}-\frac{2n-1}{4})}2^{k(\frac{2n+1}{2\widetilde{q}}-\frac{2n-1}{4})} 2^{\frac{-|j-k|}{2}(\frac{1}{2}-\max(\frac{1}{q},\frac{1}{\widetilde{q}}))}
        \|H\|_{L^{\widetilde{q}'}_{t}\mathfrak{L}^{\widetilde{q}'}_{p}}\|\widetilde{H}\|_{L^{q'}_{t} \mathfrak{L}^{q'}_{r}}
    \end{align*}
for $2\leq q,\widetilde{q}\leq6$.
By using this and \eqref{1st_cond_modify}, the first part where $j\leq k$ is now bounded as follows:
    \begin{align}\label{interpolation_value_2}
    \nonumber&\sum_{j=0}^{\infty}\Big| \sum_{k=j}^{\infty} B_{j,k} (H,\widetilde{H}) \Big| \\
    \nonumber&\lesssim\sum_{j=0}^{\infty} \sum_{k=j}^{\infty}
        2^{j(\frac{2n+1}{2q}-\frac{2n-1}{4})}2^{k(\frac{2n+1}{2\widetilde{q}}-\frac{2n-1}{4})}
        2^{\frac{-|j-k|}{2}(\frac{1}{2}-\max(\frac{1}{q},\frac{1}{\widetilde{q}}))}
        \|H\|_{L^{\widetilde{q}'}_t \mathfrak{L}^{\widetilde{q}'}_r}\|\widetilde{H}\|_{L^{q'}_t \mathfrak{L}^{q'}_r}  \\
    \nonumber&=\sum_{j=0}^{\infty}
        2^{j(\frac{2n+1}{2q}-\frac{2n-2}{4}-\frac{1}{2}\max(\frac{1}{q},\frac{1}{\widetilde{q}}))}
        \sum_{k=j}^{\infty} 2^{k(\frac{2n+1}{2\widetilde{q}}-\frac{n}{2}+ \frac{1}{2}\max(\frac{1}{q},\frac{1}{\widetilde{q}}))}
        \|H\|_{L^{\widetilde{q}'}_t \mathfrak{L}^{\widetilde{q}'}_r}
        \|\widetilde{H}\|_{L^{q'}_t \mathfrak{L}^{q'}_r}\\
    &\lesssim \sum_{j=0}^{\infty}2^{j(\frac{2n+1}{2}(\frac{1}{q}+\frac{1}{\widetilde{q}})-\frac{2n-1}{2})}
        \|H\|_{L^{\widetilde{q}'}_t \mathfrak{L}^{\widetilde{q}'}_r}\|\widetilde{H}\|_{L^{q'}_t \mathfrak{L}^{q'}_r}
    \end{align}
for $2(n+1)/n<q,\widetilde{q}\leq6$.
If one applies this bound directly for $q,\widetilde{q}$ satisfying the conditions in Proposition \ref{prop2} as in the previous case,
then one can not sum over $j$ because
$\frac{2n+1}{2}(\frac{1}{q}+\frac{1}{\widetilde{q}}) - \frac{2n-1}{2}=0$ when $\alpha = 2n/(2n-1)$.
But here we will make use of the following bilinear interpolation lemma
(see \cite{BL}, Section 3.13, Exercise 5(b))
together with \eqref{interpolation_value_2} to give
\begin{equation}\label{goal_bilinear_sum}
\sum_{j=0}^{\infty} \Big| \sum_{k=j}^{\infty} B_{j,k} (H,\widetilde{H}) \Big|
\lesssim \|H\|_{L^{\widetilde{q}'}_t \mathfrak{L}^{\widetilde{q}'}_r}
\|\widetilde{H}\|_{L^{q'}_t \mathfrak{L}^{q'}_r}.
\end{equation}

\begin{lem}\label{bilinear_lem}
For $i=0,1$, let $A_i,B_i,C_i$ be Banach spaces  and let $T$ be a bilinear operator such that
    \begin{align*}
    &T:A_0\times B_0\rightarrow C_0,\\
    &T:A_0\times B_1\rightarrow C_1,\\
    &T:A_1\times B_0\rightarrow C_1.
    \end{align*}
Then one has, for $\theta=\theta_0+\theta_1$ and $1/q+1/r\geq1$,
    $$
    T:(A_0,A_1)_{\theta_0,q}\times(B_0,B_1)_{\theta_1,r}\rightarrow(C_0,C_1)_{\theta,1}.
    $$
Here, $0<\theta_i<\theta<1$ and $1\leq q,r\leq\infty$.
\end{lem}

Indeed, let us first consider the vector-valued bilinear operator $T$ defined by
    \begin{equation*}
    T(H,\widetilde{H}) = \big\{ T_j (H,\widetilde{H}) \big\}_{j\geq0},
    \end{equation*}
where $T_j=\sum_{k=j}^{\infty} B_{j,k}$.
Then, \eqref{goal_bilinear_sum} is equivalent to
    \begin{equation}\label{bii}
    T: L^{\widetilde{q}'}_t \mathfrak{L}^{\widetilde{q}'}_r\times L^{q'}_t \mathfrak{L}^{q'}_r
     \rightarrow \ell_1^0(\mathbb{C}),
    \end{equation}
where $\ell^a_p(\mathbb{C})$, $a\in\mathbb{R}$, $1\leq p\leq\infty$, denotes
the weighted sequence space equipped with the norm
    $$
    \|\{x_j\}_{j\geq0} \|_{\ell^a_p} =
        \begin{cases}
        \big(\sum_{j\geq0}2^{jap}|x_j|^p\big)^{\frac{1}{p}}
        \quad\text{if}\quad p\neq\infty,\\
        \,\sup_{j\geq0}2^{ja}|x_j|
        \quad\text{if}\quad p=\infty.
        \end{cases}
    $$
Now, by \eqref{interpolation_value_2} we see that
    \begin{equation}\label{357}
    \| T(H,\widetilde{H}) \|_{\ell_\infty^{\beta(\widetilde{q},q)}(\mathbb{C})}
    \lesssim \|H\|_{L^{\widetilde{q}'}_t \mathfrak{L}^{\widetilde{q}'}_r}
    \|\widetilde{H}\|_{L^{q'}_t\mathfrak{L}^{q'}_r},
    \end{equation}
where $2(n+1)/n<q,\widetilde{q}\leq6$ and
$$\beta(\widetilde{q},q) = \frac{2n-1}{2} - \frac{2n+1}{2}(\frac{1}{q}+\frac{1}{\widetilde{q}}).$$
Also, for $(\widetilde{q},q)$ satisfying \eqref{accep_and_scal_cond},
we can take a sufficiently small $\epsilon>0$ such that
the ball $B((\frac{1}{\widetilde{q}},\frac{1}{q}),3\epsilon)$
with center $(\frac{1}{\widetilde{q}},\frac{1}{q})$ and radius $3\epsilon$
is contained in the region of
$(\frac{1}{\widetilde{q}},\frac{1}{q})$
given by $2(n+1)/n<q,\widetilde{q}\leq6$ (see Figure \ref{fig2}).
Now, we choose $\widetilde{q}_0, \widetilde{q}_1, q_0, q_1$ such that
    $$
    \frac{1}{\widetilde{q}_0}= \frac{1}{\widetilde{q}}-\epsilon,\quad \frac{1}{\widetilde{q}_1}= \frac{1}{\widetilde{q}}+2\epsilon, \quad
    \frac{1}{q_0}= \frac{1}{q}-\epsilon,\quad \frac{1}{q_1}= \frac{1}{q}+2\epsilon.
    $$
Then it is easy to check that
    $$
    \beta(\widetilde{q}_0,q_0) = (2n+1)\epsilon ,
    \quad \beta(\widetilde{q}_0,q_1) = - \frac{2n+1}{2}\epsilon ,
    \quad \beta(\widetilde{q}_1,q_0) = - \frac{2n+1}{2}\epsilon ,
    $$
and we get from \eqref{357} the following three bounds
    \begin{align*}
    &T: L^{\widetilde{q}_0'}_t\mathfrak{L}^{\widetilde{q}_0'}_r\times
    L^{q_0'}_t\mathfrak{L}^{q_0'}_r\rightarrow\ell^{(2n+1)\epsilon}_{\infty}(\mathbb{C}),\\
    &T: L^{\widetilde{q}_0'}_t \mathfrak{L}^{\widetilde{q}_0'}_r\times
    L^{q_1'}_t\mathfrak{L}^{q_1'}_r\rightarrow\ell^{-\frac{2n+1}{2}\epsilon}_{\infty}(\mathbb{C}),\\
    &T: L^{\widetilde{q}_1'}_t \mathfrak{L}^{\widetilde{q}_1'}_r\times
    L^{q_0'}_t\mathfrak{L}^{q_0'}_r\rightarrow\ell^{-\frac{2n+1}{2}\epsilon}_{\infty}(\mathbb{C}).
    \end{align*}
Then, by applying Lemma \ref{bilinear_lem} with $\theta_0=\theta_1=1/3$ and $q=r=2$, we get
$$T:(L^{\widetilde{q}_0'}_t\mathfrak{L}^{\widetilde{q}_0'}_r,L^{\widetilde{q}_1'}_t \mathfrak{L}^{\widetilde{q}_1'}_r)_{1/3,2}\times(L^{q_0'}_t\mathfrak{L}^{q_0'}_r, L^{q_1'}_t\mathfrak{L}^{q_1'}_r)_{1/3,2}
\rightarrow(\ell^{(2n+1)\epsilon}_{\infty}(\mathbb{C}),
\ell^{-\frac{2n+1}{2}\epsilon}_{\infty}(\mathbb{C}))_{2/3,1}.$$
Since $L_t^p\mathfrak{L}_r^p=L_{t,r}^p(r^{n-1}drdt)$,
by applying the real interpolation space identities in the following lemma,
one can easily deduce \eqref{bii} from the above boundedness.

\begin{lem}[Theorems 5.4.1 and 5.6.1 in \cite{BL}]
Let $0<\theta<1$ and $1\leq q_0,q_1,q\leq\infty$.
Then one has
$$(L^{q_0}(w_0), L^{q_1}(w_1))_{\theta,2}=L^q(w)$$
if $\frac{1}{q}=\frac{1-\theta}{q_0}+\frac{\theta}{q_1}$ and $w(x)= w_0(x)^{1-\theta}w_1(x)^{\theta}$,
and if $s=(1-\theta)s_0+\theta s_1$
$$(\ell^{s_0}_{q_0}, \ell^{s_1}_{q_1})_{\theta, q}=\ell^s_q.$$
\end{lem}

It is clear from the same argument that the second part where $j\geq k$ is bounded as follows:
$$\sum_{k=0}^{\infty} \Big| \sum_{j=k}^{\infty} B_{j,k} (H,\widetilde{H}) \Big|
\lesssim \|H\|_{L^{q'}_t \mathfrak{L}^{q'}_r}\|\widetilde{H}\|_{L^{\widetilde{q}'}_t \mathfrak{L}^{\widetilde{q}'}_r}.$$
Consequently, we have the desired estimate \eqref{2:1}.

\subsection{Proof of Lemma \ref{local_prop}}\label{subsec3.2}

It remains to prove Lemma \ref{local_prop}.
We have to prove the following estimate:
For $j,k\geq0$,
\begin{align}\label{lemm}
\nonumber\|T_j T^*_{k}&(\lambda^{n-1}H)\|_{L^{q}_{t}\mathfrak{L}^{q}_{r}}\\
&\lesssim2^{j(\frac{2n+1}{2q}-\frac{2n-1}{4})}2^{k(\frac{2n+1}{2\widetilde{q}}-\frac{2n-1}{4})} 2^{\frac{-|j-k|}{2}(\frac{1}{2}-\max(\frac{1}{q},\frac{1}{\widetilde{q}}))}
\|H\|_{L^{\widetilde{q}'}_t\mathfrak{L}^{\widetilde{q}'}_r}
\end{align}
if $2\leq q,\widetilde{q}\leq6$.
The proof is divided into the case $j,k\geq1$ and the cases where $j=0$ or $k=0$.

\

\noindent\textbf{The case $j,k \geq 1$.}
First we claim that for $2\leq q,\widetilde{q}\leq6$,
    \begin{equation}\label{inter_value_1}
    \|T_j T^{*}_{k}(\lambda^{n-1}H)\|_{L^{q}_{t}\mathfrak{L}^{q}_{r}}
    \lesssim 2^{j(\frac{2n+1}{2q}-\frac{2n-1}{4})}2^{k(\frac{2n+1}{2\widetilde{q}}-\frac{2n-1}{4})}
    \|H\|_{L^{\widetilde{q}'}_t\mathfrak{L}^{\widetilde{q}'}_{r}}.
    \end{equation}
Indeed, we note that for $j\geq1$ and $2\leq q\leq6$,
    \begin{equation}\label{ChoLee_basic_result}
    \| T_j h \|_{L^{q}_{t} \mathfrak{L}^{q}_{r}}
    \lesssim 2^{j(\frac{2n+1}{2q} - \frac{2n-1}{4})} \|h\|_2
    \end{equation}
which follows immediately from applying Proposition 3.1 in \cite{CL}
with $\varpi(\rho)=-\rho^\alpha$, $R\sim2^j$, and $p=q$.
Then by the usual $TT^\ast$ argument, it is not difficult to see that
\eqref{ChoLee_basic_result} implies
\begin{equation*}
\|T_jT^{*}_{k}H\|_{L^{q}_{t}\mathfrak{L}^{q}_{r}}
\lesssim 2^{j(\frac{2n+1}{2q}-\frac{2n-1}{4})}2^{k(\frac{2n+1}{2\widetilde{q}}-\frac{2n-1}{4})}
\|H\lambda^{-(n-1)}\|_{L^{\widetilde{q}'}_s\mathfrak{L}^{\widetilde{q}'}_{\lambda}}
\end{equation*}
for $j,k\geq1$ and $2\leq q,\widetilde{q}\leq6$.
Replacing $H$ with $\lambda^{n-1}H$, this gives \eqref{inter_value_1}.

When $|j-k|\leq1$, \eqref{lemm} follows now directly from \eqref{inter_value_1}.
So we are reduced to showing \eqref{lemm} when $|j-k|>1$.
For this we will obtain
    \begin{equation}\label{van_der_Corput_result}
    \| T_jT^{*}_{k}(\lambda^{n-1}H)\|_{L^{\infty}_{t}\mathfrak{L}^{\infty}_{r}}
    \lesssim2^{-j\frac{2n-1}{4}}2^{-k\frac{2n-1}{4}}2^{-\frac{1}{4}|j-k|}
    \|H\|_{L^{1}_t \mathfrak{L}^{1}_r}
    \end{equation}
when $|j-k|>1$.
By interpolation between the estimates \eqref{inter_value_1} and \eqref{van_der_Corput_result},
we then get \eqref{lemm} when $|j-k|>1$.
Indeed, when $\max(\frac{1}{q},\frac{1}{\widetilde{q}})=\frac{1}{q}$,
by interpolation between \eqref{inter_value_1} with $q=2$ and \eqref{van_der_Corput_result},
it is easy to check that \eqref{lemm} holds for $2\leq q\leq\infty$ and $\widetilde{q}/3\leq q\leq\widetilde{q}$. This range of $q,\widetilde{q}$ is wider than what is given by
$2\leq q,\widetilde{q}\leq6$.
When $\max(\frac{1}{q},\frac{1}{\widetilde{q}})=\frac{1}{\widetilde{q}}$,
one can similarly get \eqref{lemm} by interpolation
between \eqref{inter_value_1} with $\widetilde{q}=2$ and \eqref{van_der_Corput_result}.

Now it remains to show the estimate \eqref{van_der_Corput_result}.
From \eqref{111} and \eqref{112}, we first write
    $$
    T_j T^{*}_{k} H (r,t) =
    \iint K(r,\lambda,t-s) H(\lambda,s) d\lambda ds,
    $$
where $K(r,\lambda,t)$ is given as
    $$
    K(r,\lambda,t)
    = \frac{\chi_{[2^{j-1},2^j)}(r)}{ r^{\frac{n-2}{2}} } \frac{\chi_{[2^{k-1},2^k)}(\lambda)}{ \lambda^{\frac{n-2}{2}} }
        \int e^{it\rho^{\alpha}} J_{\frac{n-2}{2}} (r\rho) J_{\frac{n-2}{2}} (\lambda\rho) \varphi^2 (\rho) d\rho .
    $$
Then, \eqref{van_der_Corput_result} would follow from the uniform bound
    \begin{equation}\label{L^infty_kernel}
    \|K\|_{L^{\infty}_{r,\lambda,t}}
    \lesssim 2^{-j\frac{2n-1}{4}}2^{-k\frac{2n-1}{4}}2^{-\frac{1}{4}|j-k|}.
    \end{equation}
To show this bound, we will divide $K$ into four parts
based on the following estimates for Bessel functions $J_{\nu}(r)$.

\begin{lem}\label{bessel}
For $r>1$ and $Re\,\nu>-1/2$,
\begin{eqnarray*}
J_{\nu}(r)=\frac{\sqrt{2}}{\sqrt{\pi r}}\cos(r-\frac{\nu\pi}2-\frac{\pi}4)
-\frac{(\nu-\frac12)\Gamma(\nu+\frac32)}{(2\pi)^{\frac12}(r)^{\frac32}\Gamma(\nu+\frac12)}
\sin(r-\frac{\nu\pi}2-\frac{\pi}4)+E_{\nu}(r),
\end{eqnarray*}
where
\begin{equation}\label{Bessel_error_1}
|E_{\nu}(r)|\leq
C_\nu r^{-\frac52}
\end{equation}
and
\begin{equation}\label{Bessel_error_2}
|\frac{d}{dr}E_{\nu}(r)|\leq C_\nu
(r^{-\frac52}+r^{-\frac72}).
\end{equation}
\end{lem}
Assuming this lemma which will be shown in Section \ref{sec5}, we see that
\begin{equation}\label{er}
J_{\frac{n-2}{2}}(\lambda \rho)=
(c_{n}(\lambda\rho)^{-\frac{1}{2}}+c_{n}(\lambda\rho)^{-\frac{3}{2}})e^{\pm i\lambda\rho}
+E_{\frac{n-2}{2}} (\lambda\rho)
\end{equation}
and
$$J_{\frac{n-2}{2}}(r \rho)=
(c_{n}(r\rho)^{-\frac{1}{2}}+c_{n}(r\rho)^{-\frac{3}{2}})e^{\pm ir\rho}+E_{\frac{n-2}{2}}(r\rho),$$
where the letter $c_{n}$ stands for constants different at each occurrence
and depending only on $n$.
Now we write
$$J_{\frac{n-2}{2}}(\lambda\rho) J_{\frac{n-2}{2}}(r\rho)= \sum_{l=1}^4J_l(r,\lambda,\rho),$$
where
 \begin{align*}
J_1(r,\lambda,\rho)&=
(c_{n}(\lambda\rho)^{-\frac{1}{2}}+c_{n}(\lambda\rho)^{-\frac{3}{2}})
e^{\pm i\lambda\rho}(c_{n}(r\rho)^{-\frac{1}{2}}+c_{n}(r\rho)^{-\frac{3}{2}} )e^{\pm ir\rho},\\ J_2(r,\lambda,\rho)&=
(c_{n}(\lambda\rho)^{-\frac{1}{2}}+c_{n}(\lambda\rho)^{-\frac{3}{2}})
e^{\pm i\lambda\rho}E_{\frac{n-2}{2}}(r\rho),\\
J_3(r,\lambda,\rho)&=
(c_{n}(r\rho)^{-\frac{1}{2}}+c_{n}(r\rho)^{-\frac{3}{2}})e^{\pm ir\rho}E_{\frac{n-2}{2}}(\lambda\rho),\\
J_4(r,\lambda,\rho)&=E_{\frac{n-2}{2}}(\lambda\rho)E_{\frac{n-2}{2}}(r\rho).
\end{align*}
Then, $K$ is divided as $K=\sum_{l=1}^4K_l$, where
$$K_l(r,\lambda,t)=
\frac{\chi_{[2^{j-1},2^j)}(r)}{r^{\frac{n-2}{2}}}
\frac{\chi_{[2^{k-1},2^k)}(\lambda)}{\lambda^{\frac{n-2}{2}}}
\int e^{it\rho^{\alpha}}J_l(r,\lambda,\rho)\varphi^2(\rho)d\rho.$$

First, it follows easily from \eqref{Bessel_error_1} that
    $$
    \|K_4\|_{L^{\infty}_{r,\lambda,t} }
    \lesssim 2^{-j\frac{n+3}{2}}2^{-k\frac{n+3}{2}}
    \leq 2^{-j\frac{2n-1}{4}}2^{-k\frac{2n-1}{4}}2^{-\frac{1}{4}|j-k|}.
    $$

Next, we shall consider $K_1$.
Since the factors $(\lambda\rho)^{-\frac32}$ and $(r\rho)^{-\frac32}$in $J_1$
would give a better boundedness
than $(\lambda\rho)^{-\frac12}$ and $(r\rho)^{-\frac12}$, respectively,
we only need to show the bound \eqref{L^infty_kernel} for
    $$
    \widetilde{K}_1(r,\lambda,t)
    = \frac{\chi_{[2^{j-1},2^j)}(r)}{ r^{\frac{n-2}{2}} } \frac{\chi_{[2^{k-2},2^k)}(\lambda)}{ \lambda^{\frac{n-2}{2}} }(\lambda r)^{-\frac12}
        \int e^{it\rho^\alpha\pm i\lambda\rho \pm ir\rho} \rho^{-1} \varphi^2(\rho)d\rho.
    $$
Let us now decompose $\widetilde{K}_1$ as
\begin{align*}
\widetilde{K}_1(r,\lambda,t)&=
\chi_{\{\frac{2^{m(j,k)}}{8}<|t|<8\cdot2^{m(j,k)}\}}(t)\widetilde{K}_1(r,\lambda,t)\\
&+(1-\chi_{\{\frac{2^{m(j,k)}}{8}<|t|<8\cdot2^{m(j,k)}\}}(t))\widetilde{K}_1(r,\lambda,t),
\end{align*}
where $m(j,k)=\max(j,k)$.
When $\frac{2^{m(j,k)}}{8}<|t|<8\cdot2^{m(j,k)}$, by the van der Corput lemma (see \cite{St}, Chap. VIII)
it follows that
$$\Big|\int e^{it\rho^\alpha\pm i\lambda\rho \pm ir\rho} \rho^{-1} \varphi^2(\rho)d\rho\Big|\lesssim
2^{-\frac{1}{2}m(j,k)}.$$
Hence, we get
\begin{align*}
\big\|\chi_{\{\frac{2^{m(j,k)}}{8}<|t|<8\cdot2^{m(j,k)}\}}(t)\widetilde{K}_1\big\|_{L^{\infty}_{r,\lambda,t}}
&\lesssim2^{-k\frac{n-1}{2}}2^{-j\frac{n-1}{2}}2^{-\frac{1}{2}m(j,k)}\\
&=2^{-j\frac{2n-1}{4}}2^{-k\frac{2n-1}{4}}2^{-\frac{1}{4}|j-k|}.
\end{align*}
For the second part where $|t|>8\cdot2^{m(j,k)}$ or $|t|<\frac{2^{m(j,k)}}{8}$,
we first note that
\begin{align*}
&\int \Big(1-\frac{\alpha(\alpha-1)t\rho^{\alpha-2}}{i(\pm\lambda \pm r +\alpha t\rho^{\alpha-1})^2}\Big)
e^{\pm i\lambda\rho \pm ir\rho +it\rho^\alpha} \rho^{-1}\varphi^2(\rho)d\rho\\
&\qquad\qquad=\Big[\frac{1}{i(\pm\lambda \pm r +\alpha t\rho^{\alpha-1})}e^{\pm i\lambda\rho \pm ir\rho +it\rho^\alpha}
\rho^{-1} \varphi^2 (\rho)  \Big]_{\rho=1/2}^{\rho=2}\\
&\qquad\qquad-\int \frac{1}{i(\pm\lambda \pm r+ \alpha t\rho^{\alpha-1})} e^{\pm i\lambda\rho \pm ir\rho +it\rho^\alpha} \frac{d}{d\rho} \Big( \rho^{-1}\varphi^2 (\rho) \Big) d\rho
\end{align*}
by integration by parts.
Since we are handling the case where $|j-k|>1$,
we see that $|\pm\lambda \pm r+\alpha t\rho^{\alpha-1}|\gtrsim 2^{m(j,k)}$
when $|t|>8\cdot2^{m(j,k)}$ or $|t|<\frac{2^{m(j,k)}}{8}$.
Hence, using this, we get
    \begin{align*}
    \Big|\int e^{\pm i\lambda\rho \pm ir\rho +it\rho^\alpha} \rho^{-1}\varphi^2(\rho)d\rho\Big|
    &\leq\int \Big| \frac{\alpha(\alpha-1)t\rho^{\alpha-2}}{(\pm\lambda \pm r+ \alpha t\rho^{\alpha-1})^2} \rho^{-1}\varphi^2 (\rho) \Big| d\rho \\
    &\quad +\Big|\Big[\frac{1}{(\pm\lambda \pm r+\alpha t\rho^{\alpha-1})}
    \rho^{-1}\varphi^2 (\rho) \Big]_{\rho=1/2}^{\rho=2} \Big|\\
    &\quad +\int\Big|\frac{1}{(\pm\lambda \pm r+\alpha t\rho^{\alpha-1})}\frac{d}{d\rho}
    \Big(\rho^{-1}\varphi^2 (\rho) \Big) \Big| d\rho\\
    &\lesssim2^{-m(j,k)}.
    \end{align*}
This implies
\begin{align*}
\|(1-\chi_{\{\frac{2^{m(j,k)}}{8}<|t|<8\cdot2^{m(j,k)}\}}(t))\widetilde{K}_1\|_{L^{\infty}_{r,\lambda,t} }
&\lesssim 2^{-k\frac{n-1}{2}} 2^{-j\frac{n-1}{2}} 2^{-m(j,k)} \\
&\leq 2^{-j\frac{2n-1}{4}}2^{-k\frac{2n-1}{4}}2^{-\frac{1}{4}|j-k|}.
\end{align*}

It remains to bound $K_2$ and $K_3$.
We shall show the bound \eqref{L^infty_kernel} only for $K_2$
because the same type of argument used for $K_2$ works clearly on $K_3$.
Since the factor $(\lambda\rho)^{-\frac32}$ in $J_2$ would give a better boundedness
than $(\lambda\rho)^{-\frac12}$, we only need to show the bound \eqref{L^infty_kernel} for
    $$
    \widetilde{K}_2(r,\lambda,t)
    = \frac{\chi_{[2^{j-1},2^j)}(r)}{ r^{\frac{n-2}{2}} } \frac{\chi_{[2^{k-1},2^k)}(\lambda)}{ \lambda^{\frac{n-2}{2}} }
        \int e^{it\rho^\alpha\pm i\lambda\rho}E_{\frac{n-2}{2}}(r\rho) (\lambda\rho)^{-\frac12}\varphi^2(\rho)d\rho.
    $$
Let us now decompose $\widetilde{K}_2$ as
\begin{align*}
\widetilde{K}_2(r,\lambda,t)&=
\chi_{\{\frac{2^{m(j,k)}}{8}<|t|<8\cdot2^{m(j,k)}\}}(t)\widetilde{K}_2(r,\lambda,t)\\
&+(1-\chi_{\{\frac{2^{m(j,k)}}{8}<|t|<8\cdot2^{m(j,k)}\}}(t))\widetilde{K}_2(r,\lambda,t),
\end{align*}
where $m(j,k)=\max(j,k)$.
When $\frac{2^{m(j,k)}}{8}<|t|<8\cdot2^{m(j,k)}$, by the van der Corput lemma as before,
it follows that
\begin{align*}
|\chi_{\{\frac{2^{m(j,k)}}{8}<|t|<8\cdot2^{m(j,k)}\}}(t)&\widetilde{K}_2(r,\lambda,t)|\\
\lesssim&\frac{\chi_{[2^{k-1},2^k)}(\lambda)\chi_{[2^{j-1},2^j)}(r)}
{\lambda^{\frac{n-2}{2}}r^{\frac{n-2}{2}}}\lambda^{-\frac12}2^{-\frac{1}{2}m(j,k)}\\
\times&\sup_{\rho\in[1/2,2]}\Big\{E_{\frac{n-2}{2}}(r\rho)\rho^{-\frac12}\varphi^2(\rho),
\frac{d}{d\rho}\big(E_{\frac{n-2}{2}}(r\rho)\rho^{-\frac12}\varphi^2(\rho)\big)\Big\}.
\end{align*}
By \eqref{Bessel_error_1} and \eqref{Bessel_error_2} in Lemma \ref{bessel},
we see
    \begin{equation}\label{Bessel_error_con}
    | E_{\frac{n-2}{2}} (r\rho) | \lesssim r^{-\frac{5}{2}}
    \quad\text{and}\quad
    | \frac{d}{d\rho} E_{\frac{n-2}{2}} (r\rho) | \lesssim r^{-\frac{3}{2}}
    \end{equation}
for $1/2 \leq \rho \leq 2$.
Thus, we get
\begin{align*}
\big\|\chi_{\{\frac{2^{m(j,k)}}{8}<|t|<8\cdot2^{m(j,k)}\}}(t)\widetilde{K}_2\big\|_{L^{\infty}_{r,\lambda,t}}
&\lesssim 2^{-k\frac{n-1}{2}}2^{-j\frac{n-2}{2}}2^{-\frac{1}{2}m(j,k)}2^{-j\frac{3}{2}}\\
&\leq 2^{-j\frac{2n-1}{4}}2^{-k\frac{2n-1}{4}}2^{-\frac{1}{4}|j-k|}.
\end{align*}
For the second part where $|t|>8\cdot2^{m(j,k)}$ or $|t|<\frac{2^{m(j,k)}}{8}$,
we will use the following trivial bound when $m(j,k)=j$ and $r\sim2^j$:
    \begin{align*}
    \Big|\int e^{it\rho^\alpha\pm i\lambda\rho} E_{\frac{n-2}{2}} (r\rho)\rho^{-\frac12} \varphi^2 (\rho) d\rho \Big|
    &\leq\int \big| E_{\frac{n-2}{2}} (r\rho)\rho^{-\frac12} \varphi^2 (\rho) \big| d\rho \\
    &\lesssim 2^{-\frac{5}{2}j} = 2^{-m(j,k)} 2^{-\frac32j}
    \end{align*}
which follows from \eqref{Bessel_error_con}.
On the other hand, when $m(j,k)=k$ and $r\sim2^j$, we will also show
$$
\Big|\int e^{it\rho^\alpha\pm i\lambda\rho}E_{\frac{n-2}{2}}(r\rho)\rho^{-\frac12}\varphi^2(\rho)d\rho\Big|
\lesssim 2^{-m(j,k)} 2^{-\frac{3}{2}j}.
$$
Indeed, by integration by parts we see that
\begin{align*}
&\int \Big(1-\frac{\alpha(\alpha-1)t\rho^{\alpha-2}}{i(\pm\lambda+\alpha t\rho^{\alpha-1})^2}\Big)
e^{\pm i\lambda\rho+it\rho^\alpha}E_{\frac{n-2}{2}}(r\rho)\rho^{-\frac12}\varphi^2(\rho)d\rho\\
&\qquad\qquad=\Big[\frac{1}{i(\pm\lambda+\alpha t\rho^{\alpha-1})}e^{\pm i\lambda\rho+it\rho^\alpha}
E_{\frac{n-2}{2}} (r\rho)\rho^{-\frac12} \varphi^2 (\rho)  \Big]_{\rho=1/2}^{\rho=2}\\
&\qquad\qquad-\int \frac{1}{i(\pm\lambda+ \alpha t\rho^{\alpha-1})} e^{\pm i\lambda\rho+it\rho^\alpha} \frac{d}{d\rho} \Big( E_{\frac{n-2}{2}} (r\rho) \rho^{-\frac12}\varphi^2 (\rho) \Big) d\rho.
\end{align*}
Since $m(j,k)=k$, one can easily check that $|\pm\lambda+\alpha t\rho^{\alpha-1}|\gtrsim 2^{m(j,k)}$
when $|t|>8\cdot2^{m(j,k)}$ or $|t|<\frac{2^{m(j,k)}}{8}$.
Hence, using this and \eqref{Bessel_error_con}, we get
\begin{align*}
\Big|\int e^{\pm i\lambda\rho+it\rho^\alpha}&E_{\frac{n-2}{2}}(r\rho)\rho^{-\frac12}\varphi^2(\rho)d\rho\Big|\\
&\leq\int \Big| \frac{\alpha(\alpha-1)t\rho^{\alpha-2}}{(\pm\lambda+ \alpha t\rho^{\alpha-1})^2} E_{\frac{n-2}{2}} (r\rho) \rho^{-\frac12}\varphi^2 (\rho) \Big| d\rho \\
&+\Big|\Big[\frac{1}{(\pm\lambda+\alpha t\rho^{\alpha-1})}E_{\frac{n-2}{2}}(r\rho)
\rho^{-\frac12}\varphi^2 (\rho) \Big]_{\rho=1/2}^{\rho=2} \Big|\\
&+\int\Big|\frac{1}{(\pm\lambda+\alpha t\rho^{\alpha-1})}\frac{d}{d\rho}
\Big(E_{\frac{n-2}{2}} (r\rho) \rho^{-\frac12}\varphi^2 (\rho) \Big) \Big| d\rho\\
&\lesssim2^{-m(j,k)}2^{-\frac{3}{2}j}
\end{align*}
as desired.
Consequently, if $r\sim2^j$
$$\Big|\int e^{it\rho^\alpha\pm i\lambda\rho}
E_{\frac{n-2}{2}} (r\rho)\rho^{-\frac12} \varphi^2 (\rho) d\rho \Big|
\lesssim2^{-m(j,k)} 2^{-\frac32j}$$
when $|t|>8\cdot2^{m(j,k)}$ or $|t|<\frac{2^{m(j,k)}}{8}$.
This implies
\begin{align*}
\|(1-\chi_{\{\frac{2^{m(j,k)}}{8}<|t|<8\cdot2^{m(j,k)}\}}(t))\widetilde{K}_2\|_{L^{\infty}_{r,\lambda,t} }
&\lesssim 2^{-k\frac{n-1}{2}} 2^{-j\frac{n-2}{2}} 2^{-m(j,k)} 2^{-\frac32j}\\
&\leq 2^{-j\frac{2n-1}{4}}2^{-k\frac{2n-1}{4}}2^{-\frac{1}{4}|j-k|}.
\end{align*}

\

\noindent\textbf{The cases where $j=0$ or $k=0$.}
Now we consider the following cases where $j=0$ or $k=0$ in \eqref{lemm}:
\begin{equation}\label{remain_claim_1}
    \| T_0 T^*_{k} (\lambda^{n-1}H) \|_{L^{q}_{t}\mathfrak{L}^{q}_{r}}
    \lesssim 2^{k(\frac{2n+1}{2\widetilde{q}}-\frac{2n-1}{4})}2^{\frac{-k}{2} \big( \frac{1}{2}-\max(\frac{1}{q},\frac{1}{\widetilde{q}}) \big)}
    \|H\|_{L^{\widetilde{q}'}_{t}\mathfrak{L}^{\widetilde{q}'}_{r}},
    \end{equation}
  \begin{equation}\label{remain_claim_2}
    \| T_j T^{*}_{0} (\lambda^{n-1}H) \|_{L^{q}_{t}\mathfrak{L}^{q}_{r}}
    \lesssim 2^{j(\frac{2n+1}{2q}-\frac{2n-1}{4})} 2^{\frac{-j}{2} \big( \frac{1}{2}-\max(\frac{1}{q},\frac{1}{\widetilde{q}}) \big)}
    \|H\|_{L^{\widetilde{q}'}_{t}\mathfrak{L}^{\widetilde{q}'}_{r}},
   \end{equation}
and
\begin{equation}\label{T0_T*0_case}
    \big\| T_0 T^{*}_0 (\lambda^{n-1}H) \big\|_{L^{q}_{t} \mathfrak{L}^{q}_{r}}
    \lesssim \|H\|_{L^{\widetilde{q}'}_t \mathfrak{L}^{\widetilde{q}'}_r},
    \end{equation}
where $j,k\geq1$ and $2\leq q,\widetilde{q}\leq6$.

Since the second estimate \eqref{remain_claim_2} follows easily from the first one
using the dual characterisation of $L^p$ spaces and a property of adjoint operators,
we only show \eqref{remain_claim_1} and \eqref{T0_T*0_case} repeating the previous argument.
But here we use the following estimates (see \cite{G}, p. 426) for Bessel functions instead of Lemma \ref{bessel}:
For $0\leq r<1$ and $Re\,\nu>-1/2$,
\begin{equation}\label{Besse}
|J_{\nu} (r)| \leq C_{\nu} r^{\nu}
\quad\text{and}\quad
\big|\frac{d}{dr}J_{\nu} (r)\big| \leq C_{\nu} r^{\nu-1}.
\end{equation}

First we shall show \eqref{T0_T*0_case}.
Recall that
    $$
    T_0h(t,r)=\chi_{(0,1)}(r) r^{-\frac{n-2}2} \int_0^\infty e^{it\rho^{\alpha}}J_{\frac{n-2}2}(r\rho)\varphi(\rho)h(\rho) d\rho.
    $$
Then, by changing variables $\rho=\rho^{\alpha}$, we see that
    $$
    \int_0^\infty e^{it\rho^{\alpha}}J_{\frac{n-2}2}(r\rho)\varphi(\rho)h(\rho) d\rho
    =\alpha^{-1}\int_0^\infty e^{it\rho}J_{\frac{n-2}2}(r\rho^{1/{\alpha}})\varphi(\rho^{1/{\alpha}})h(\rho^{1/{\alpha}}) \rho^{1/{\alpha} - 1}d\rho.
    $$
Thus, using Plancherel's theorem in $t$ and \eqref{Besse}, we get
\begin{align}\label{T0_L2-L2_estimate}
\nonumber\| T_0 h \|_{L_t^2 \mathcal{L}_r^2}
&= \| T_0 h \|_{\mathcal{L}_r^2 L_t^2} \\
\nonumber&= C\Big\| \chi_{(0,1]}(r) r^{-\frac{n-2}2} J_{\frac{n-2}2}(r\rho^{1/{\alpha}})
\varphi(\rho^{1/{\alpha}})h(\rho^{1/{\alpha}})\rho^{1/{\alpha} - 1} \Big\|_{\mathcal{L}_r^2 L_\rho^2} \\
\nonumber&= C\Big( \int_0^1 r^{-(n-2)} \int_{1/2}^2 |J_{\frac{n-2}2}(r\rho)|^2|h(\rho)|^2 \rho^{1 - {\alpha}} d\rho r^{n-1}dr \Big)^{1/2}\\
\nonumber&\lesssim\Big( \int_{1/2}^2 |h(\rho)|^2 \rho^{1-{\alpha}} \int_0^1 C r(r\rho)^{n-2}dr d\rho \Big)^{1/2}\\
\nonumber&\lesssim\Big(\int_{1/2}^2 |h(\rho)|^2 d\rho\Big)^{1/2}\\
&=\|h\|_{L^2}.
\end{align}
Also, by H\"older's inequality,
    \begin{equation*}
    \| T_0 h \|_{L_t^{\infty} \mathcal{L}_r^{\infty}}\lesssim \| h \|_{2}.
    \end{equation*}
By interpolation between this and \eqref{T0_L2-L2_estimate}, we obtain
    \begin{equation}\label{T0_estimate}
    \| T_0 h \|_{L^{q}_{t}\mathfrak{L}^{q}_{r}} \lesssim \|h\|_2
    \end{equation}
for $2 \leq q \leq \infty$.
Then, by the usual $TT^\ast$ argument as before, this implies
    \begin{equation*}
    \| T_0 T^*_0 (\lambda^{n-1}H) \|_{L^{q}_{t}\mathfrak{L}^{q}_{r}}
    \lesssim \|H\|_{L^{\widetilde{q}'}_{t}\mathfrak{L}^{\widetilde{q}'}_{r}}
    \end{equation*}
for $2 \leq q, \widetilde{q} \leq \infty$.

Now we turn to \eqref{remain_claim_1}.
First, by using \eqref{T0_estimate} and the dual estimate of \eqref{ChoLee_basic_result}, we see that
    \begin{align*}
    \| T_0 T^*_{k} (\lambda^{n-1}H) \|_{L^{q}_{t}\mathfrak{L}^{q}_{r}}
    &\lesssim \| T^*_{k} (\lambda^{n-1}H) \|_2 \\
    &\lesssim 2^{k(\frac{2n+1}{2\widetilde{q}} - \frac{2n-1}{4})} \|H\|_{L^{\widetilde{q}'}_{t}\mathfrak{L}^{\widetilde{q}'}_{r}}
    \end{align*}
for $2 \leq q \leq \infty$ and $2 \leq \widetilde{q} \leq 6$.
Then, \eqref{remain_claim_1} would follow from interpolation between this and the following estimate
as before (see the paragraph below \eqref{van_der_Corput_result}):
     \begin{equation}\label{T_0T^*_k estimate}
    \| T_0 T^{*}_{k} (\lambda^{n-1}H) \|_{L^{\infty}_{t} \mathfrak{L}^{\infty}_{r}}
    \lesssim 2^{-k\frac{n}{2}} \|H\|_{L^{1}_t \mathfrak{L}^{1}_r}.
    \end{equation}
Now we are reduced to showing \eqref{T_0T^*_k estimate}.
From \eqref{110} and \eqref{112}, we first write
    $$
    T_0 T^{*}_{k} H (r,t) =
    \iint K(r,\lambda,t-s) H(\lambda,s) d\lambda ds,
    $$
where
    $$
    K(r,\lambda,t)
    = \frac{\chi_{(0,1)}(r)}{ r^{\frac{n-2}{2}} } \frac{\chi_{[2^{k-1},2^k)}(\lambda)}{ \lambda^{\frac{n-2}{2}}}
        \int e^{it\rho^{\alpha}} J_{\frac{n-2}{2}} (r\rho) J_{\frac{n-2}{2}} (\lambda\rho) \varphi^2 (\rho) d\rho .
    $$
Then, we only need to show that
   \begin{equation}\label{inft}
    \|K(r,\lambda,t)\|_{L^\infty_{r,\lambda,t}} \lesssim 2^{-k\frac{n}{2}}.
    \end{equation}
Recall from \eqref{er} that
\begin{equation}\label{jj}
J_{\frac{n-2}{2}}(\lambda \rho)=
(c_{n}(\lambda\rho)^{-\frac{1}{2}}+c_{n}(\lambda\rho)^{-\frac{3}{2}})e^{\pm i\lambda\rho}
+E_{\frac{n-2}{2}} (\lambda\rho).
\end{equation}
By \eqref{Bessel_error_con} and \eqref{Besse},
the part of $K$ coming from $E_{\frac{n-2}{2}}(\lambda\rho)$ in $\eqref{jj}$ is bounded as follows:
\begin{equation*}
\Big|\frac{\chi_{(0,1)}(r)}{r^{\frac{n-2}{2}}}\frac{\chi_{[2^{k-1},2^k)}(\lambda)}{\lambda^{\frac{n-2}{2}}}
\int e^{it\rho^{\alpha}}J_{\frac{n-2}{2}}(r\rho)E_{\frac{n-2}{2}}(\lambda\rho)\varphi^2(\rho) d\rho\Big|\lesssim2^{-k\frac{n+3}{2}}\leq2^{-k\frac{n}{2}}.
\end{equation*}
Now we may consider only the part of $K$ coming from $(\lambda\rho)^{-\frac{1}{2}}$,
because the factor $(\lambda\rho)^{-\frac32}$ in $\eqref{jj}$ would give a better boundedness
than $(\lambda\rho)^{-\frac12}$.
Namely, we have to show the bound \eqref{inft} for
$$\widetilde{K}(r,\lambda,t)
=\frac{\chi_{(0,1)}(r)}{r^{\frac{n-2}{2}}}\frac{\chi_{[2^{k-1},2^k)}(\lambda)}{\lambda^{\frac{n-2}{2}}}
\int e^{it\rho^{\alpha}\pm i\lambda\rho}
J_{\frac{n-2}{2}}(r\rho)(\lambda\rho)^{-\frac12}\varphi^2(\rho)d\rho.$$
Let us now decompose $\widetilde{K}$ as
$$\widetilde{K}(r,\lambda,t)=
\chi_{\{\frac{2^k}{8}<|t|<8\cdot2^k\}}(t)\widetilde{K}(r,\lambda,t)
+(1-\chi_{\{\frac{2^k}{8}<|t|<8\cdot2^k\}}(t))\widetilde{K}(r,\lambda,t).$$
When $\frac{2^k}{8}<|t|<8\cdot2^k$, by the van der Corput lemma as before,
it follows that
\begin{align*}
|\chi_{\{\frac{2^k}{8}<|t|<8\cdot2^k\}}(t)&\widetilde{K}(r,\lambda,t)|\\
\lesssim&\frac{\chi_{[2^{k-1},2^k)}(\lambda)\chi_{(0,1)}(r)}
{\lambda^{\frac{n-2}{2}}r^{\frac{n-2}{2}}}\lambda^{-\frac12}2^{-\frac{1}{2}k}\\
\times&\sup_{\rho\in[1/2,2]}\Big\{J_{\frac{n-2}{2}}(r\rho)\rho^{-\frac12}\varphi^2(\rho),
\frac{d}{d\rho}\big(J_{\frac{n-2}{2}}(r\rho)\rho^{-\frac12}\varphi^2(\rho)\big)\Big\}.
\end{align*}
By \eqref{Besse}, we see
\begin{equation*}
|J_{\frac{n-2}{2}} (r\rho)| \leq C r^{\frac{n-2}{2}}
\quad\text{and}\quad \big|\frac{d}{d\rho}J_{\frac{n-2}{2}} (r\rho)\big| \leq C r^{\frac{n-2}{2}}
\end{equation*}
for $1/2 \leq \rho \leq 2$.
Thus, we get
$$
\big\|\chi_{\{\frac{2^k}{8}<|t|<8\cdot2^k\}}(t)\widetilde{K}\big\|_{L^{\infty}_{r,\lambda,t}}
\lesssim 2^{-k\frac{n}{2}}.
$$
On the other hand, by integration by parts we see that
\begin{align*}
&\int\Big(1-\frac{\alpha(\alpha-1)t\rho^{\alpha-2}}{i(\pm\lambda+\alpha t\rho^{\alpha-1})^2}\Big)
e^{\pm i\lambda\rho+it\rho^\alpha}J_{\frac{n-2}{2}}(r\rho)\rho^{-\frac12}\varphi^2(\rho)d\rho\\
&\qquad\qquad=\Big[\frac{1}{i(\pm\lambda+\alpha t\rho^{\alpha-1})}e^{\pm i\lambda\rho+it\rho^\alpha}
J_{\frac{n-2}{2}}(r\rho)\rho^{-\frac12}\varphi^2(\rho)\Big]_{\rho=1/2}^{\rho=2}\\
&\qquad\qquad-\int\frac{1}{i(\pm\lambda+\alpha t\rho^{\alpha-1})}e^{\pm i\lambda\rho+it\rho^\alpha} \frac{d}{d\rho}\Big(J_{\frac{n-2}{2}}(r\rho)\rho^{-\frac12}\varphi^2(\rho)\Big)d\rho.
\end{align*}
One can also easily check that $|\pm\lambda+\alpha t\rho^{\alpha-1}|\gtrsim 2^k$
when $|t|>8\cdot2^k$ or $|t|<\frac{2^k}{8}$.
Hence, using this and \eqref{Besse}, we get
\begin{align*}
\Big|\int e^{\pm i\lambda\rho+it\rho^\alpha}&J_{\frac{n-2}{2}}(r\rho)\rho^{-\frac12}\varphi^2(\rho)d\rho\Big|\\
&\leq\int \Big| \frac{\alpha(\alpha-1)t\rho^{\alpha-2}}{(\pm\lambda+ \alpha t\rho^{\alpha-1})^2} J_{\frac{n-2}{2}} (r\rho) \rho^{-\frac12}\varphi^2 (\rho) \Big| d\rho \\
&+\Big|\Big[\frac{1}{(\pm\lambda+\alpha t\rho^{\alpha-1})}J_{\frac{n-2}{2}}(r\rho)
\rho^{-\frac12}\varphi^2 (\rho) \Big]_{\rho=1/2}^{\rho=2} \Big|\\
&+\int\Big|\frac{1}{(\pm\lambda+\alpha t\rho^{\alpha-1})}\frac{d}{d\rho}
\Big(J_{\frac{n-2}{2}} (r\rho) \rho^{-\frac12}\varphi^2 (\rho) \Big) \Big| d\rho\\
&\lesssim2^{-k}r^{\frac{n-2}{2}}
\end{align*}
when $|t|>8\cdot2^k$ or $|t|<\frac{2^k}{8}$.
This implies
$$\|(1-\chi_{\{\frac{2^k}{8}<|t|<8\cdot2^k\}}(t))\widetilde{K}\|_{L^{\infty}_{r,\lambda,t} }
\lesssim 2^{-k\frac{n-2}{2}}2^{-k}=2^{-k\frac n2}.$$

\section{Sharpness of Theorem \ref{prop}}\label{sec4}

In this section we discuss the sharpness of Theorem \ref{prop}.
We will show that \eqref{nece2} is a necessary condition for \eqref{inhomo} (see Remark \ref{rem}).
If \eqref{nece2} is valid with a pair $(q,p)$ on the left and a pair $(\widetilde{q},\widetilde{p})$
on the right, then it must be also valid when one switches the roles of $(q,p)$ and $(\widetilde{q},\widetilde{p})$.
By this duality relation, we only need to show the first condition $n/p+1/q<n/2$ in \eqref{nece2}.

Let $\phi$ be a smooth cut-off function supported in the interval $[1,2]$.
Let us now define $F(y,s)$ by
$$\widehat{F(\cdot,s)}(\xi)=\chi_{(0,1)}(s)\phi(|\xi|).$$
Then one can easily see that $\|F\|_{L^{\widetilde{q}'}_{s}L^{\widetilde{p}'}_{y}}\lesssim1$,
and
\begin{align*}
\Big\|\int_{0}^{t}e^{i(t-s)|\nabla|^{\alpha}}&F(\cdot,s)ds\Big\|_{L^{q}_tL^{p}_x}\\
&\geq\Big\|\int_{0}^{1}\int e^{ix\cdot\xi+i(t-s)|\xi|^{\alpha}}\phi(|\xi|)d\xi ds \Big\|_{L^{q}_{t}((N,\infty);L^{p}_{x}(\frac{t}{2}<|x|<t))}\\
&=\Big\|\int e^{ix\cdot\xi+it|\xi|^{\alpha}}\Big(\frac{e^{-i|\xi|^{\alpha}}-1}{-i|\xi|^{\alpha}} \Big)\phi(|\xi|)d\xi\Big\|_{L^{q}_{t}((N,\infty);L^{p}_{x}(\frac{t}{2}<|x|<t))}
\end{align*}
by taking integration with respect to $s$ as
$$\int_{0}^{1} e^{-is|\xi|^{\alpha}} ds = \frac{e^{-i|\xi|^\alpha}-1}{-i|\xi|^\alpha}.$$

Now we recall from \cite{St} (see p. 344 there) that
$$I(\lambda):=\int e^{i\lambda\psi(\xi)}\omega(\xi)d\xi\sim\lambda^{-\frac{n}2} \sum_{j=0}^{\infty}a_j\lambda^{-j}$$
where $a_0\neq0$, $\psi$ has a nondegenerate critical point at some point $\xi_0$
(i.e., $\nabla \psi(\xi_0)=0$ and the matrix $\big[\frac{\partial^2\psi}{\partial\xi_i\partial\xi_j}\big](\xi_0)$ is invertible),
and $\omega$ is supported in a sufficiently small neighborhood of $\xi_0$.
Applying this with $\psi(\xi)=\frac1t x\cdot\xi+|\xi|^{\alpha}$ and
$$\omega(\xi)=\Big(\frac{-e^{i|\xi|^{\alpha}}-1}{-i|\xi|^{\alpha}}\Big)\phi(|\xi|),$$
we get
$$\Big|\int e^{ix\cdot\xi+it|\xi|^{\alpha}}\Big(\frac{-e^{i|\xi|^{\alpha}}-1}{-i|\xi|^{\alpha}}\Big)\phi(|\xi|)d\xi
\Big|\gtrsim t^{-\frac{n}{2}}$$
for sufficiently large $t$.
Thus, if $N$ is sufficiently large,
    \begin{align*}
    \Big\| \int e^{ix\cdot\xi + it|\xi|^{\alpha}} \Big(\frac{-e^{i|\xi|^{\alpha}}-1}{-i|\xi|^{\alpha}}\Big)&\phi(|\xi|)
    d\xi\Big\|_{L^{q}_{t}((N,\infty ; L^{p}_{x} (\frac{t}{2} < |x| < t))}\\
    &\gtrsim \bigg( \int_{N}^{\infty} t^{-\frac{n}{2}q} \Big( \int_{\frac{t}{2} < |x| < t } dx \Big)^{\frac{q}{p}} dt \bigg)^{\frac{1}{q}} \\
    &\sim \Big( \int_{N}^{\infty} t^{nq ( \frac{1}{p}-\frac{1}{2} ) } dt \Big)^{\frac{1}{q}}.
    \end{align*}
Consequently, if \eqref{inhomo} holds,
$$\Big(\int_{N}^{\infty}t^{nq(1/p-1/2)}dt\Big)^{\frac{1}{q}}\lesssim1.$$
But, this is not possible as $N\rightarrow\infty$ unless $nq(1/p-1/2)<-1$
which is equivalent to the condition $n/p+1/q<n/2$.

\section{Appendix}\label{sec5}

Here we shall provide a proof of Lemma \ref{bessel} for estimates of Bessel functions $J_{\nu}(r)$.
It is based on easy but quite tedious calculations.

First, we recall from \cite{G} (see p. 430 there) that for $r>1$ and $Re\,\nu>-1/2$,
$$J_{\nu}(r)=\frac{(r/2)^{\nu}}{\Gamma(\nu+\frac12)\Gamma(\frac12)}
\Big(ie^{-ir}\int_0^{\infty}e^{-rt}(t^2+2it)^{\nu-\frac12}dt  -ie^{ir}\int_0^{\infty}e^{-rt}(t^2-2it)^{\nu-\frac12}dt\Big),$$
where $\Gamma$ is the gamma function given by
\begin{equation}\label{ga}
\Gamma(z) = \int_0^{\infty} x^{z-1} e^{-x} dx,\quad Re\,z>0.
\end{equation}
Then, using the following identities
$$ie^{-ir}(t^2+2it)^{\nu-\frac12} = (2t)^{\nu-\frac12}e^{-i(r-\frac{\nu\pi}2-\frac{\pi}4)}(1-\frac{it}2)^{\nu-\frac12}$$
and
$$-ie^{ir}(t^2-2it)^{\nu-\frac12} = (2t)^{\nu-\frac12}e^{i(r-\frac{\nu\pi}2-\frac{\pi}4)}(1+\frac{it}2)^{\nu-\frac12}$$
together with
\begin{equation}\label{firsttay2}
(1-\frac{it}2)^{\nu-\frac12} = 1 - (\nu-\frac12)\frac{it}2 + R_{\nu}(t)
\end{equation}
and
\begin{equation}\label{firsttay1}
(1+\frac{it}2)^{\nu-\frac12} = 1 + (\nu-\frac12)\frac{it}2 + \widetilde{R}_{\nu}(t),
\end{equation}
one can rewrite
\begin{align*}
J_{\nu}(r)=\frac{(2\pi)^{-\frac12}r^{\nu}}{\Gamma(\nu+\frac12)}
&\bigg(e^{-i(r-\frac{\nu\pi}2-\frac{\pi}4)}\int_0^{\infty}e^{-rt}t^{\nu-\frac12}
\Big(1 -(\nu-\frac12)\frac{it}{2}+R_{\nu}(t)\Big)dt\\
&+e^{i(r-\frac{\nu\pi}2-\frac{\pi}4)}\int_0^{\infty}e^{-rt}t^{\nu-\frac12}
\Big(1 + (\nu-\frac12)\frac{it}{2}+\widetilde{R}_{\nu}(t)\Big)dt\bigg).
\end{align*}
Here, $R_{\nu}(t)$ and $\widetilde{R}_{\nu}(t)$ are the remainder terms
in Taylor series \eqref{firsttay2} and \eqref{firsttay1}, respectively,
which are given by
$$R_{\nu}(t)=(\nu-\frac12)(\nu-\frac32)(\frac{it}2)^2(1-\frac{it_{*}}{2})^{\nu-\frac52}$$
and
$$\widetilde{R}_{\nu}(t)=(\nu-\frac12)(\nu-\frac32)(\frac{it}2)^2(1+\frac{it^{*}}{2})^{\nu-\frac52}$$
for some $t_{*}$ and $t^{*}$ with $0<t_{*},t^{*}<t$.

Now we decompose $J_{\nu}(r)$ into three parts as $J_{\nu}(r)=I+II+III$, where
$$I = \frac{(2\pi)^{-\frac12}r^{\nu}}{\Gamma(\nu+\frac12)}
(e^{-i(r-\frac{\nu\pi}2-\frac{\pi}4)}+e^{i(r-\frac{\nu\pi}2-\frac{\pi}4)})
\int_0^{\infty}e^{-rt}t^{\nu-\frac12}dt,$$
$$II=\frac{(\nu-\frac12)(2\pi)^{-\frac12}r^{\nu}}{2\Gamma(\nu+\frac12)}
(-ie^{-i(r-\frac{\nu\pi}2-\frac{\pi}4)}+ie^{i(r-\frac{\nu\pi}2-\frac{\pi}4)})
\int_0^{\infty}e^{-rt}t^{\nu+\frac12}dt$$
and
$$III=\frac{(2\pi)^{-\frac12}r^{\nu}}{\Gamma(\nu+\frac12)}
\Big(\int_0^{\infty}\frac{e^{-rt}t^{\nu-\frac12}R_{\nu}(t)}{e^{i(r-\frac{\nu\pi}2-\frac{\pi}4)}}dt  +\int_0^{\infty}\frac{e^{-rt}t^{\nu-\frac12}
\widetilde{R}_{\nu}(t)}{e^{-i(r-\frac{\nu\pi}2-\frac{\pi}4)}}dt\Big).$$
Then, from the definition \eqref{ga} of $\Gamma$, we easily see that
$$I=\frac{\sqrt{2}}{\sqrt{\pi r}}\cos(r-\frac{\nu\pi}2-\frac{\pi}4)$$
and
$$II=-\frac{(\nu-\frac12)\Gamma(\nu+\frac32)}{(2\pi)^{\frac12}r^{\frac32}
\Gamma(\nu+\frac12)}\sin(r-\frac{\nu\pi}2-\frac{\pi}4).$$

Now, the lemma is proved by taking $E_{\nu}(r)=III$.
Indeed, to show \eqref{Bessel_error_1} and \eqref{Bessel_error_2},
we first note that when $Re\,\nu-5/2\leq0$,
\begin{equation}\label{remainder1}
|R_{\nu}(t)|,\, |\widetilde{R}_{\nu}(t)|\leq C_\nu|(\nu-\frac12)(\nu-\frac32)|t^2,
\end{equation}
and when $Re\,\nu-5/2>0$,
\begin{equation}\label{remainder2}
|R_{\nu}(t)|,\, |\widetilde{R}_{\nu}(t)| \leq C_\nu|(\nu-\frac12)(\nu-\frac32)|t^{2}\max(t^{\nu-\frac52},1).
\end{equation}
Hence, when $Re\,\nu-5/2\leq0$,
by using \eqref{remainder1} and \eqref{ga}, it follows that
$$|III| \leq C_\nu|\nu-\frac12||\nu-\frac32|\frac{\Gamma(\nu+\frac52)}
{\Gamma(\nu+\frac12)}r^{-\frac52}.$$
On the other hand, when $Re\,\nu-5/2>0$,
by using \eqref{remainder2} and \eqref{ga},
we see that
\begin{align*}
|III|&\leq C_\nu|\nu-\frac12||\nu-\frac32|\frac{\max(\Gamma(\nu+\frac52)r^{-\frac52},\Gamma(2\nu)r^{-\nu})}
{\Gamma(\nu+\frac12)}\\
&\leq C_\nu|\nu-\frac12||\nu-\frac32|\frac{\max(\Gamma(\nu+\frac52),\Gamma(2\nu))}{\Gamma(\nu+\frac12)}r^{-\frac52}.
\end{align*}
Consequently, we get
\begin{equation*}
|E_{\nu}(r)|\leq C_\nu r^{-\frac52},
\end{equation*}
and \eqref{Bessel_error_2} is similarly shown by
differentiating $E_{\nu}(r)$ and using \eqref{remainder1} and \eqref{remainder2}.


\end{document}